\input amstex.tex
\documentstyle{amsppt}
\input pictex.tex
\magnification=\magstep1
\pagewidth{16.1truecm}
\pageheight{23.9truecm}
\nologo
\NoRunningHeads

\def\sing{\operatorname{sing}}
\def\logdisc{\operatorname{logdisc}}
\def\totallogdisc{\operatorname{totallogdisc}} 
\def\Var{\operatorname{Var}}
\def\supp{\operatorname{supp}}
\def\dim{\operatorname{dim}}
\def\and{\operatorname{and}}
\def\deg{\operatorname{deg}}
\def\nnc{\operatorname{nnc}}
\def\Diff{\operatorname{Diff}}

\refstyle{A}
\widestnumber\key{AKMW}
\topmatter
\title
Stringy zeta functions for $\Bbb Q$--Gorenstein varieties
\endtitle
\author
Willem Veys  \bigskip
\endauthor
\address K.U.Leuven, Departement Wiskunde, Celestijnenlaan 200B,
         B--3001 Leuven, Belgium  \endaddress
\email wim.veys\@wis.kuleuven.ac.be  \newline
 http://www.wis.kuleuven.ac.be/algebra/veys.htm
\endemail
\abstract
The stringy Euler number and stringy E--function are interesting
invariants of log terminal singularities, introduced by Batyrev. He used them to formulate a topological mirror symmetry test for pairs of certain Calabi--Yau varieties, and to show a version of the McKay correspondence. It is a natural question whether one can extend these invariants beyond the log terminal case.
Assuming the Minimal Model Program, we introduce very general stringy
invariants, associated to \lq almost all\rq\ singularities, more precisely to all singularities which are not strictly log canonical.
They specialize to the invariants of Batyrev when the singularity is log terminal. For example the simplest form of our stringy zeta function is in general a rational function in one variable, but it is just a constant (Batyrev's stringy Euler number) in the log terminal case.
\endabstract
\endtopmatter
\bigskip
\document

\heading
Introduction
\endheading
\bigskip
\noindent
{\bf 0.1.} The stringy Euler number and stringy $E$--function are interesting singularity invariants introduced by Batyrev.  In [Ba1] they are associated to log terminal complex algebraic varieties $X$, and in [Ba2] more generally to Kawamata log terminal pairs $(X,B)$.  Batyrev used them to formulate a topological mirror symmetry test for pairs of certain Calabi--Yau varieties, and to show a version of the McKay correspondence. They are also the subject of remarkable conjectures [Ba1].  We recall their definition.

Let $X$ be a normal complex variety and $B = \sum_i b_i B_i$ a  $\Bbb Q$--divisor on $X$, where the $B_i$ are distinct and irreducible, such that $K_X + B$ is $\Bbb Q$--Cartier.  Here $K_\centerdot$  denotes the canonical divisor. (In particular when $B=0$ this means that $X$ is $\Bbb Q$--Gorenstein.)  

For a birational morphism 
$\pi : Y \rightarrow X$ from a normal variety $Y$, let $E_i, i \in T$, denote the irreducible divisors in the union of $\pi^{-1}(\supp B)$ and the exceptional locus of $\pi$.  The {\it log discrepancies} $a_i$ of $E_i, i \in T$, with respect to $(X,B)$ are given by
$$K_Y = \pi^\ast (K_X + B) + \sum_{i \in T} (a_i - 1)E_i.$$
The pair $(X,B)$ is Kawamata log terminal (klt) precisely when for a log resolution $\pi : Y \rightarrow X$ of the pair $(X,B)$ we have that $a_i > 0$ for all $i \in T$.  This condition does not depend on the chosen resolution.  Remark that if $E_i$ is the strict transform of a component $B_i$, then $a_i = 1 - b_i$; hence all $b_i < 1$ for a klt pair.

Denote also $E^\circ_I :=(\cap_{i \in I} E_i) \setminus (\cup_{\ell \not\in I} E_\ell)$ for $I \subset T$.  So $Y$ is the disjoint union of the $E^\circ_I, I \subset T$. Then the {\it stringy Euler number} of the klt pair $(X,B)$ is
$$e(X,B) := \sum_{I \subset T} \chi (E^\circ_I) \prod_{i \in I} \frac{1}{a_i} \in \Bbb Q ,$$   
where $\chi(\cdot)$ denotes the topological Euler characteristic.  A finer invariant is the {\it stringy $E$--function}
$$E(X,B) := \sum_{I \subset T} H(E^\circ_I)  \prod_{i \in I} \frac{uv-1}{(uv)^{a_i}-1}$$
of $(X,B)$, where $H(E^\circ_I) \in \Bbb Z[u,v]$ is the Hodge polynomial of $E^\circ_I$, see (1.8).  (When $B=0$ we just write $e(X)$ and $E(X)$. For smooth $X$ we have $e(X)=\chi(X)$ and $E(X)=H(X)$.)  The proof of Batyrev that these invariants do not depend on the chosen resolution uses the idea of motivic integration, initiated by Kontsevich [Kon] and developed by Denef and Loeser [DL2][DL3][DL4].  Another proof is possible using weak factorization [AKMW][W{\l}].

We refer to [Da] and [DR] for concrete formulas for $E(X)$ and $e(X)$ in some interesting cases, and [BM, \S2] for a comparison with the string--theoretic $E$--polynomial of [BD]. 
See also (1.8.1) for a relation between $E(X)$ and the singular elliptic genus of $X$.
\bigskip
\noindent
{\bf 0.2.} It is a natural question whether one can generalize these invariants beyond the log terminal case.  For surfaces $X$ we realized this in [Ve3], see (3.5).

In arbitrary dimension we obtained in [Ve2, \S3] invariants, given by the same formulas, for pairs $(X,B)$, where $X$ is any $\Bbb Q$--Gorenstein variety and $B$ an effective $\Bbb Q$--Cartier divisor on $X$ whose support contains $X_{\sing}$, or, more generally, contains the locus of log canonical singularities of $X$.  Of course, working with $B = 0$, this is not more general than in (0.1).   
\bigskip
\noindent
{\bf 0.3.} In complete generality, we think that the following two questions are the most natural to ask.  Let $(X,B = \sum_i b_iB_i)$ be any pair as in (0.1) with all $b_i < 1$.

(I) Suppose that there exists a log resolution $Y$ of $(X,B)$ on which all log discrepancies with respect to $(X,B)$ are nonzero.  Then one can state the same formulas as in (0.1) using $Y$.  If $Y^\prime$ is another such log resolution, do the formulas using $Y$ and $Y^\prime$ give the same expression ?  A positive answer would yield stringy invariants for pairs $(X,B)$ admitting at least one such resolution.  The natural approach to this question using weak factorization encounters an annoying difficulty, just as in [BL2] for the elliptic genus, see (5.6).

(II) Do there exist invariants, associated to {\it any} $\Bbb Q$--Gorenstein variety $X$, and more generally to any pair $(X,B)$, that specialize to the stringy invariants of (0.1) in case $X$ is log terminal and $(X,B)$ is klt, respectively ?

\bigskip
\noindent
{\bf 0.4.} In this paper we just mention some partial results concerning question (I) in (3.4), (4.4) and \S 5.  Our main purpose is to attack question (II) assuming the Minimal Model Program (MMP).  More precisely we will assume the relative log MMP and associate {\it stringy zeta functions} on different levels as in (0.1) to $\Bbb Q$--Gorenstein varieties $X$ and to pairs $(X,B)$.
\bigskip
\noindent
{\bf 0.5.} For example on the roughest level of Euler characteristics we associate an invariant $z(X;s) \in \Bbb Q(s)$ to `almost all' $\Bbb Q$--Gorenstein varieties $X$, such that for log terminal $X$ we have that $z(X;s) = e(X)$ and is thus a constant.

We present this more in detail.  Let $p : X^m \rightarrow X$ be a relative log minimal model of $X$, see (1.6).  We consider the generic case where all log discrepancies with respect to $X$ of exceptional divisors on $X^{m}$ are negative.  In (3.7) we explain that this condition is indeed `very generic'; it is equivalent to $X$ being either log terminal (then there are simply no exceptional divisors on $X^{m}$), or not log canonical and without strictly log canonical singularities.  On $X^{m}$ we use the divisor 
$D := (K_{X^m} + F) - p^\ast K_X$, where $F = \sum_{i \in T^m}  F_i$ is the (reduced) exceptional divisor of $p$.  By definition of the log discrepancies with respect to $X$, we have that $D = \sum_{i \in T^m} a_i F_i$, where all $a_i < 0$ by assumption.

Take now a log resolution $h : Y \rightarrow X^m$ of the pair $(X^m,F)$ such that $\pi = p \circ h : Y \rightarrow X$ is a log resolution of $X$.  We use for $\pi$ the notation of (0.1); in particular $E_i, i \in T^m \subset T$, is the strict transform by $h$ of $F_i$.  We associate to $X$ the {\it stringy zeta function}
$$z(X;s) := \sum_{I \subset T} \chi (E^\circ_I) \prod_{i \in I} \frac{1}{\nu_i + sN_i} ,$$
where
$$\align K_Y & = h^\ast (K_{X^m} + F) + \sum_{i \in T} (\nu_i - 1)E_i \qquad \text{ and }\\
 h^\ast D & = \sum_{i \in T} N_i E_i; \endalign$$
thus in particular the $\nu_i$ are the log discrepancies of the $E_i$ with respect to the pair $(X^m, F)$. 

The expression on the right hand side at least makes sense : for the exceptional components $E_i, i \in T \setminus T^m$, of $h$ we have by definition of a relative log minimal model that either $\nu_i > 0$, or $\nu_i = 0$ and $N_i < 0$, and for the $E_i$, $i \in T^m$, we have that $\nu_i = 0$ but $N_i = a_i < 0$.  We will show that this expression does not depend on the choices of $X^m$ and $Y$.

The fact that some $\nu_i$ are necessarily zero, is a technical problem that is in some sense forced by nature: the pair $(X^m,F)$ is not klt but only divisorial log terminal (dlt), see (1.4) and (1.6).

\smallskip
It is easy to see that for $i \in T$ the log discrepancy $a_i = \nu_i + N_i$.  So if there exists a log resolution $Y$ of $X$ with all $a_i \ne 0$, factorizing through some $X^m$ (by a morphism), then $z(X;1) \in \Bbb Q$ is independent of such $Y$ and is given by the same formula as in (0.1).

In particular when $X$ is log terminal, we have that  $T_m=\emptyset$, hence $D=0$, all $N_i=0$ and $\nu_i=a_i$, and $z(X;s)$ is indeed just Batyrev's $e(X)$.

\bigskip
\noindent
{\bf 0.6.}  This generalization is consistent with our definition of the stringy Euler number for {\sl surfaces} $X$ in [Ve3].  There we associate $e(X) \in \Bbb Q$ to {\sl any} normal surface $X$ without strictly log canonical singularities, obtaining another generalization of Batyrev's expression in (0.1).

That approach consists essentially in the following : define explicitly the contribution to $e(X)$ of exceptional curves $E_\ell$ with $a_\ell =0$, and use the same contribution as in (0.1) for the strata $E^\circ_I$ which are disjoint with those $E_\ell$, see (3.5).
It turns out that the relevant contribution to define is for $E_\ell \cong \Bbb P^1$ with $a_\ell=0$, intersecting two other components, say $E_1$ and $E_2$, with $a_1\neq 0 \neq a_2$. And the \lq right\rq\ contribution of this $E_\ell$ is then $\frac{-E_\ell^2}{a_1a_2}$.
(See [NN] for a topological interpretation of this generalized $e(X)$ in a special case.)

We show more precisely
in Proposition 3.5.4 that our $e(X)$ for such normal surfaces equals precisely $z(X;1)$, confirming the naturality of both approaches.  We stress that it is a priori not clear that this evaluation $z(X;1)$ makes sense because some $a_i, i \in T \setminus T^m$, really can be zero !

\bigskip
\noindent
{\bf 0.7.} Let $X$ be any $\Bbb Q$-Gorenstein variety without strictly log canonical singularities. Considering (0.6) and the last paragraph of (0.5) one could propose
$$
\lim_{s\to 1} z(X;s) \in \Bbb Q \cup \{\infty\}
$$
as a definition for a generalized stringy Euler number $e(X)$. 
The real question here is whether (or when) this limit exists in $\Bbb Q$, as in the surface case. Also, is it then possible to define explicit contributions of the $E_\ell$ with $a_\ell =0$ to such an $e(X)$ ?
We do not know the answer but we think this is worth to investigate.

We present in (3.6) a concrete \lq positive example\rq\ of a threedimensional $X$ with an exceptional surface $E_\ell$ with log discrepancy $a_\ell =0$ in some log resolution, and such that $\lim_{s\to 1} z(X;s) \in \Bbb Q$.

\bigskip
\noindent
{\bf 0.8.}  Everything in (0.5) also works for pairs $(X,B)$.  Moreover we will associate similar stringy zeta functions to {\it any} $\Bbb Q$--Gorenstein $X$ or pair $(X,B = \sum_i b_i B_i)$ with all $b_i < 1$.  In fact for $d \in \Bbb Q$ with $0 \leq d < 1$, we introduce analogously $z_d(X;s)$ and $z_d((X,B);s)$ in $\Bbb Q(s)$ in terms of a relative $d$--minimal model of $X$ and $(X,B)$, respectively; see (1.6) for this notion.  (The `usual' relative log minimal model is a relative 1--minimal model.) An advantage of this maybe \lq less natural\rq\ object is that for $d<1$ a $d$--minimal model is klt. 

Also we will proceed on the most general level of the Grothendieck ring of algebraic varieties, i.e. we consider the {\it universal Euler characteristic\/}, which specializes to the level of Hodge polynomials and to the above presented level of Euler characteristics.  On this general level weak factorization yields a priori finer results than motivic integration, see (2.8).
For instance birationally equivalent (smooth, complete) Calabi--Yau varieties determine the same element in a {\sl localization} of the Grothendieck ring, and not just in a completion of it as described in [DL4, 4.4.2].

At this point we admit not to be aware of potential applications of our stringy zeta functions in the sense of the nice applications that Batyrev produced.

\bigskip
\noindent
{\bf 0.9.} On the level of Hodge polynomials we prove a functional equation for our stringy zeta functions, which generalizes the Poincar\'e duality result of Batyrev [Ba1, Theorem 3.7].
\bigskip
\noindent
{\bf 0.10.} We work over the base field of complex numbers.  However, everything could be generalized to an arbitrary base field of characteristic zero assuming the MMP over that field.

In $\S1$ we recall some birational geometry, in particular the notions of (relative) log minimal model and log canonical model and their $d$--variants for $d<1$, and also the weak factorization theorem. As a preparation for our stringy zeta functions we introduce in $\S2$ zeta funcions associated to arbitrary $\Bbb Q$--Cartier divisors on klt pairs, and also to certain $\Bbb Q$--Cartier divisors on dlt pairs. Their definition is in terms of a log resolution, and weak factorization is used to prove independency of the chosen resolution. Here we also verify the functional equation.

In $\S3$ we associate a stringy zeta function to a pair $(X,B)$; more precisely this will be the zeta function of $\S2$ associated to the \lq log discrepancy divisor\rq\ on a relative log minimal model of $(X,B)$. Our restriction in $\S2$ on the allowed $\Bbb Q$--divisor puts a condition on the allowed $(X,B)$. We verify that this at first sight technical condition just means that $(X,B)$ has no strictly log canonical singularities.  In this section we also prove the consistency with our previous generalizations on normal surfaces, and we present the \lq positive example\rq\ of (0.7).

The variants of these stringy zeta functions in terms of $d$--minimal models ($d<1$) are introduced in $\S4$; in fact here we can use also $d$--canonical models. We compute them explicitly for the strictly log canonical surface singularities. Finally in $\S5$ we present briefly a partial answer to question (I), restricting to log resolutions of the variety $X$ that factorize through the blowing--up of $X$ in its singular locus.

\bigskip
\bigskip
\head
{1. Birational geometry}
\endhead
\bigskip
\noindent As general references for (1.1) to (1.6) we mention [KM], [KMM] and [Kol].

\bigskip
\noindent
{\bf 1.1.} Let $X$ be an irreducible complex algebraic variety, i.e. an integral separated scheme of finite type over Spec $\Bbb C$, and $B$ a $\Bbb Q$--divisor on $X$ (we allow $B=0$).

A log resolution of the variety $X$ is a proper birational morphism $\pi : Y \rightarrow X$ from a smooth $Y$ such that the exceptional locus of $\pi$ is a (simple) normal crossings divisor.

More generally, a log resolution of the pair $(X,B)$ is a proper birational morphism $\pi : Y \rightarrow X$ from a smooth $Y$ such that the exceptional locus of $\pi$ is a divisor, and its union with $\pi^{-1}(\supp B)$ is a (simple) normal crossings divisor.

\bigskip
\noindent
{\bf 1.2.} Let moreover $X$ be normal.
A (Weil) $\Bbb Q$--divisor $D$ on $X$ is called $\Bbb Q$--Cartier if some positive integer multiple of $D$ is Cartier.  And $X$ is called $\Bbb Q$--factorial if every Weil divisor on $X$ is $\Bbb Q$--Cartier.

Let $p:X\to S$ be a proper morphism. A $\Bbb Q$--divisor $D$ on $X$ is called $p$--nef if the intersection number $D\cdot C \geq 0$ for all irreducible curves $C$ on $X$ for which $p(C)$ is a point.

The variety $X$ has a well defined linear equivalence class $K_X$ of canonical (Weil) divisors.  One says that $X$ is Gorenstein if $K_X$ is Cartier, and $\Bbb Q$--Gorenstein if $K_X$ is $\Bbb Q$--Cartier.

\bigskip
\noindent
{\bf 1.3.}  For a $\Bbb Q$--Gorenstein $X$, let $\pi : Y \rightarrow X$ be a birational morphism from a normal variety $Y$, and denote by $E_i, i \in T$, the irreducible divisors in the exceptional locus of $\pi$.  Then in the expression
$$K_Y = \pi^\ast K_X + \sum_{i \in T} (a_i - 1)E_i$$
the rational number $a_i, i \in T$, is called the {\it log discrepancy of $E_i$ with respect to $X$} (the number $a_i - 1$ is called the {\it discrepancy}).  The {\it log discrepancy} of $X$, denoted $\logdisc (X)$, is the infimum of all $a_i$, where $E_i$ runs through all irreducible exceptional divisors of all such $Y \rightarrow X$.  Either $\logdisc (X) = -\infty$, or $0 \leq \logdisc (X) \leq 2$ [KM, Corollary 2.31].
\bigskip
\noindent
{\bf 1.3.1. Definition.}  One says that $X$ is {\it terminal, canonical, log terminal and log canonical\/} if the log discrepancy of $X$ is $> 1, \geq 1, > 0$ and $\geq 0$, respectively.

In each of these cases it is sufficient to check that the $a_i, i \in T$, for one fixed log resolution $\pi : Y \rightarrow X$ of $X$ satisfy the required inequality.
\bigskip
\noindent
{\bf 1.4.} The analogous relevant notions for pairs are more subtle.  Let $X$ be a normal variety and $0 \ne B = \sum_i b_i B_i$ a $\Bbb Q$--divisor on $X$, where the $B_i$ are distinct and irreducible, such that $K_X + B$ is $\Bbb Q$--Cartier.  Again we take a birational morphism $\pi : Y \rightarrow X$ from a normal variety $Y$, and now we consider the expression
$$K_Y = \pi^\ast (K_X + B) + \sum_{i \in T} (a_i - 1)E_i . \tag $*$ $$
Here the $E_i, i \in T$, are necessarily the irreducible exceptional divisors $E_i, i \in T_e$, of $\pi$, and the strict transforms $E_i, i \in T_s$, of the $B_i$.  So $T = T_e \amalg T_s$.  Again $a_i, i \in T$, is called the {\it log discrepancy of $E_i$ with respect to $(X,B)$}.  In particular $a_i = 1 - b_i$ for $i \in T_s$.  Alternatively one can write ($*$) in the form
$$K_Y + \sum_{i \in T_s} b_i E_i + \sum_{i \in T_e} E_i = \pi^\ast (K_X + B) + \sum_{i \in T_e} a_i E_i,$$
which reflects more the comparison between the pairs $(X,B)$ and $(Y, \sum_{i \in T_s} b_i E_i + \sum_{i \in T_e} E_i$), and the naturality of the $a_i$ (versus the $a_i - 1$).  But then the log discrepancies of the strict transforms of the $B_i$ do not appear automatically.

The {\it log discrepancy} and {\it total log discrepancy} of $(X,B)$, denoted $\logdisc(X,B)$ and  $\totallogdisc (X,B)$, is the infimum of all $a_i$, where $E_i$ runs through all irreducible exceptional divisors of all $Y \rightarrow X$, and through these divisors {\it and} the strict transforms of irreducible divisors on $X$, respectively.  Either $\logdisc (X,B) = -\infty$ or $0 \leq \totallogdisc (X,B) \leq \logdisc (X,B) \leq 2$ [KM, Corollary 2.31].
\bigskip
\noindent
{\bf 1.4.1.}  Here the relevant special cases are the generalizations of log terminal and log canonical in (1.3.1). However, the `right notion' of log terminality for pairs is not clear; several ones are important in the MMP.  We mention the two notions that will be used in this paper.
\bigskip
\noindent {\bf Definition.}  (i) One says that $(X,B)$ is {\it Kawamata log terminal} (klt) if $\totallogdisc$ $(X,B) > 0$, or, equivalently, $\logdisc (X,B) > 0$ and all $b_i < 1$.

(ii) One says that $(X,B)$ is {\it log canonical} (lc) if logdisc $(X,B) \geq 0$.  (This implies that all $b_i \leq 1$ and hence that also $\totallogdisc (X,B) \geq 0$.) 

\noindent In these two cases it is again sufficient that the $a_i, i \in T$, for one fixed log resolution $\pi : Y \rightarrow X$ of $(X,B)$ satisfy $a_i > 0$ and $a_i \geq 0$, respectively.

(iii) We now suppose that all $b_i$ satisfy $0 \leq b_i \leq 1$.  One says that $(X,B)$ is {\it divisorial log terminal} (dlt) if there exists a closed $Z \subsetneqq X$ such that 

(1) $X \setminus Z$ is smooth and $B|_{X \setminus Z}$ is a normal crossings divisor, and

(2) if $\pi : Y \rightarrow X$ is a birational morphism and $E_i \subset Y$ is an irreducible divisor satisfying $\pi (E_i) \subset Z$, then $a_i > 0$.
\bigskip
\noindent
{\sl Remark.}  We may assume in (iii) that $B|_{X \setminus Z}$ is reduced, i.e. that $\cup_{b_i < 1} B_i \subset Z$, and furthermore that then $Z$ is the smallest closed subset of $X$ satisfying (1).
\bigskip
\noindent
Equivalently, $(X,B)$ is dlt if and only if there {\it exists} a log resolution $\pi : Y \rightarrow X$ of $(X,B)$ such that $a_i > 0$ for all irreducible exceptional divisors of 
$\pi$ [Sz].
\bigskip
\noindent
{\bf 1.4.2.} {\sl  Remark.}  The subtle differences between the log terminality notions are caused by the coefficients  $b_i =  1$ in $B$.  If all $b_i$ satisfy $0 \leq b_i < 1$, then $(X,B)$ is klt if and only if it is dlt [KM, Proposition 2.41].
\bigskip
\noindent
{\bf 1.4.3.} We call $P \in X$ a {\it strictly lc} singularity of $(X,B)$ if there exists a neighbourhood $U$ of $P$ such that $(U,B|_U)$ is log canonical, but there exists no neighbourhood $V$ of $P$ such that $(V,B|_V)$ is klt.  (Here we also consider $B=0$.)     
\bigskip
\noindent
{\bf 1.5.} In dimension 2 the notion of log discrepancy makes sense for any normal surface $X$ by Mumford's concept of the pull back of a Weil divisor [Mu].  In particular all notions in (1.3) and (1.4) exist for arbitrary normal surfaces $X$ and arbitrary $\Bbb Q$--divisors $B$ on $X$.
\bigskip
\noindent
{\bf 1.6.} Let $X$ be a normal variety and $B = \sum_i b_i B_i$ a $\Bbb Q$--divisor on $X$, where the $B_i$ are distinct and irreducible, and all $b_i$ satisfy $0 \leq b_i < 1$.
\bigskip
\noindent
{\bf 1.6.1. Definition.}  (1) A {\it (relative) log minimal model} of $(X,B)$ is a proper birational morphism
$$p : X^m \rightarrow X$$
such that, denoting by $F$ the reduced exceptional divisor of $p$ and by $B^m$ the strict transform of $B$ in $X^m$,

(i) $X^m$ is $\Bbb Q$--factorial,

(ii) $(X^m,B^m + F)$ is dlt, and

(iii) $K_{X^m} + B^m + F$ is $p$--nef.

\noindent
By analogy with the next notion, we call this object also a {\it (relative) 1--minimal model}.

On the other hand fix any $d \in \Bbb Q$ with $0 \leq d < 1$.  A {\it (relative) $d$--minimal model} of $(X,B)$ is a proper birational morphism
$$p : X^m_d \rightarrow X$$
such that, with analogous $F$ and $B^m$,

(i) $X_d^m$ is $\Bbb Q$--factorial,

(ii) $\logdisc (X^m_d,B^m + dF) > 1 - d$ (in particular the pair $(X^m_d, B^m + dF$) is klt), and

(iii) $K_{X^m_d} + B^m + dF$ is $p$--nef.
\bigskip
\noindent
(2) Fix $d \in \Bbb Q$ with $0 \leq d \leq 1$.  A {\it (relative) $d$--canonical model} of $(X,B)$ is a proper birational morphism
$$q : X^c_d \rightarrow X$$
such that, denoting by $F^\prime$ the reduced exceptional divisor of $q$ and by $B^c$ the strict transform of $B$ in $X^c_d$,

(i) $\logdisc (X^c_d,B^c + dF^\prime) \geq 1 - d$, and

(ii) $K_{X^c_d} + B^c + dF^\prime$ is $q$--ample.

\noindent
(For $d=1$ one rather uses the term {\it log canonical model}.)  
\bigskip
\noindent
{\bf 1.6.2.} {\sl Remarks.}  (i) The existence of these objects is essentially equivalent to the (relative, log)  Minimal Model Program.  In particular they are proved to exist in dimension $\leq 3$.  In fact one constructs them by applying the (relative, log) MMP, starting from a log resolution of $(X,B)$.  In [KM] they are called minimal and canonical model {\it of\/} this resolution {\it over} $X$.  

(ii) When $B=0$ and $d=0$ we recover the usual relative minimal and canonical model of $X$.
\bigskip
\proclaim{1.6.3. Properties} (i) For $d < 1$ two different $d$--minimal models of $(X,B)$ are isomorphic in codimension one; 
in particular $(X,B)$ has a unique $d$--minimal model when $X$ is a surface.
(This is not true for the `usual' case $d=1$.)

(ii) For any $d$ a $d$--canonical model of $(X,B)$ is unique.

(iii) Any $d$--minimal model $p : X^m_d \rightarrow X$ factors through $q : X^c_d \rightarrow X$.  
\endproclaim

\noindent
See [Kol, Theorem 6.16] and the proofs in [KM, 3.8]. (For (i) one easily verifies that two different $d$--minimal models $(d < 1)$ are both, in the terminology of [KM], minimal models of a common log resolution, and then they are isomorphic in codimension one by [KM, Theorem 3.52 (2)].)
\bigskip
\noindent
{\bf 1.7.}  We will use the weak factorization theorem of [AKMW] and [W{\l}], which is in fact valid for varieties over any field of characteristic zero.  We state it here in the form that we need.  This is more general than the statement in [AKMW] and [W{\l}], but is implicit in these papers; see Remark 1.7.2.

First we recall that a birational map $\phi : Y -\! \rightarrow Y^\prime$ is said to be proper if the projections to $Y$ and $Y^\prime$ of the graph of $\phi$ are proper.  (This reduces to the usual notion if $\phi$ is a morphism.)
\bigskip
\noindent
\proclaim{1.7.1. Theorem} (1) Let $\phi : Y -\! \rightarrow Y^\prime$ be a proper birational map between smooth irreducible varieties, and let $U \subset Y$ be an open set where $\phi$ is an isomorphism.  Then $\phi$ can be factored as follows into a sequence of blow--ups and blow--downs with smooth centres disjoint from $U$.

There exist smooth irreducible varieties $Y_1, \dots , Y_{\ell -1}$ and a sequence of birational maps
$$Y = Y_0 - \overset \phi_1 \to \rightarrow Y_1 - \overset \phi_2 \to \rightarrow \cdots -\! \overset \phi_{i-1} \to \rightarrow Y_{i-1} - \overset \phi_i \to \rightarrow Y_i - \! \overset \phi_{i+1} \to \rightarrow \cdots - \! \overset \phi_{\ell - 1} \to \rightarrow Y_{\ell - 1} -  \overset \phi_\ell \to \rightarrow Y_\ell = Y^\prime$$
where $\phi = \phi_\ell \circ \phi_{\ell - 1} \circ \cdots \circ \phi_2 \circ \phi_1$, such that each $\phi_i$ is an isomorphism over $U$ (we identify $U$ with an open in the $Y_i$), and for $i = 1, \dots , \ell$ either $\phi_i : Y_{i-1} -\! \rightarrow Y_i$ or $\phi^{-1}_i : Y_i -\! \rightarrow Y_{i-1}$ is the blowing--up at a smooth centre disjoint from $U$, and is thus a morphism.

(1$^\prime$) There is an index $i_0$ such that for all $i \leq i_0$ the map $Y_i \rightarrow Y$ is a morphism, and for $i \geq i_0$ the map $Y_i \rightarrow Y^\prime$ is a morphism.  Moreover these morphisms are all projective.

(2) If $Y \setminus U$ and $Y^\prime \setminus U$ are normal crossings divisors, then the factorization above can be chosen such that the inverse images of these divisors under $Y_i \rightarrow Y$ or $Y_i \rightarrow Y^\prime$ are also normal crossings divisors, and such that the centres of blowing--up of the $\phi_i$ or $\phi^{-1}_i$ have normal crossings with these divisors.

(3) If $Y$ and $Y^\prime$ are varieties over a base variety $S$ and $\phi$ is a map of $S$--varieties, then the factorization above is a factorization over $S$.
\endproclaim
\bigskip
\noindent
{\bf 1.7.2.} {\sl Remark.}  (i) In [AKMW] and [W{\l}] the theorem is stated for a birational map $\phi$ between complete $Y$ and $Y^\prime$; the generalization to proper birational maps between not necessarily complete $Y$ and $Y^\prime$ is mentioned by Bonavero [Bo].

(ii) In [AKMW, Theorem 0.3.1] the first claim of (2) is not explicitly stated, but can be read off from the proof (see [AKMW, 5.9 and 5.10]).

(iii) The relative statement (3) follows from (1').
\bigskip
\noindent
{\bf 1.8.} (i) The Grothendieck ring $K_0 (\Var_{\Bbb C})$ of complex algebraic varieties is the free abelian group generated by the symbols $[V]$, where $V$ is a variety, subject to the relations $[V] = [V^\prime]$, if $V$ is isomorphic to $V^\prime$, and $[V] = [V \setminus W]+ [W]$, if $W$ is closed in $V$.  Its ring structure is given by $[V] \cdot [W] := [V \times W]$.  (See [Bi] for alternative descriptions of $K_0 (\Var_{\Bbb C})$ and [Po] for the recent proof that it is not a domain.)

We abbreviate $L := [\Bbb A^1]$. For the sequel we need to extend $K_0(\Var_{\Bbb C})$ with fractional powers of $L$ and to localize with respect to elements of the form $L^q$ and $L^q - 1$.  So formally we first introduce the quotient ring
$$K_0(\Var_{\Bbb C})[L^{\Bbb Q_{>0}}] := \frac{K_0(\Var_{\Bbb C})[T_i]_{i \in \Bbb Z_{ > 0}}}{(T^i_i - L, T^k_i - T^\ell_j)_{\Sb 
i,j,k,\ell \in \Bbb Z_{ > 0}  \\ i \ell = jk \endSb}} \, .$$ 
This indeed means that we add $L^{1/i} := \overline T_i$, and of course we require that $(L^{1/i})^k = (L^{1/j})^\ell$ if $\frac k i = \frac \ell j$.

Consider then the localization of this ring with respect to the elements $L^q, q \in \Bbb Q_{>0},$ and $L^q - 1, q \in \Bbb Q \setminus \{0\}$; we denote by
$\Cal R$  its subring generated by $K_0 (\Var_{\Bbb C})[L^{\Bbb Q_{>0}}]$ and the elements $\frac{L-1}{L^q-1}, q \in \Bbb Q \setminus \{ 0 \}$.
\medskip
(ii) For a variety $V$ we denote by $h^{p,q}(H^i_c (V, \Bbb C))$ the rank of the $(p,q)$--Hodge component in the mixed Hodge structure of the $i$th cohomology group with compact support of $V$, and we put $e^{p,q}(V) := \sum_{i \geq 0} (-1)^i h^{p,q} (H^i_c(V, \Bbb C))$.  The {\it Hodge polynomial} of $V$ is
$$H(V) = H(V;u,v) := \sum_{p,q} e^{p,q} (V) u^p v^q \in \Bbb Z[u,v].$$
Precisely by the defining relations of $K_0 (\Var_{\Bbb C})$ there is a well defined ring homomorphism $H : K_0 (\Var_{\Bbb C}) \rightarrow \Bbb Z[u,v]$, determined by $[V] \mapsto H(V)$.  It induces a ring homomorphism $H$ from $\Cal R$ to the `rational functions in $u,v$ with fractional powers'.  

\medskip
(iii) The topological Euler characteristic $\chi (V)$ of a variety $V$, i.e. the alternating sum of the ranks of its Betti or de Rham cohomology groups, satisfies $\chi(V) = H(V;1,1)$ and we obtain a ring homomorphism $\chi : K_0 (\Var_{\Bbb C}) \rightarrow \Bbb Z$ determined by $[V] \mapsto \chi(V)$.  Since $\chi(L) = 1$, it induces a ring homomorphism $\chi : \Cal R \rightarrow \Bbb Q$ by declaring $\chi ( \frac{L-1}{L^q-1}) = \frac 1 q$.

\smallskip
\noindent See e.g. [DL2], [DL3] and [Ve2] for similar constructions.

\bigskip
\demo{1.8.1. Note} (1) When $X$ is projective and smooth its Hodge polynomial $H(X)$ incorporates besides $\chi(X)$ also other classical numerical invariants. We have that $H(X;-y,1) \in \Bbb Z[y]$ equals Hirzebruch's $\chi_y$--genus of $X$ [Hi], which specializes to the holomorphic Euler characteristic and signature of $X$.

Another generalization of this $\chi_y$--genus is the two variable elliptic genus of $X$ (see [BL1]); this elliptic genus of $X$ and $H(X)$ are really different generalizations in the sense that neither of both can be derived from the other one.

\smallskip
(2) When $X$ is projective and log terminal this pattern generalizes : then Borisov and Libgober defined in [BL2] a {\sl singular elliptic genus\/} of $X$, such that one of its specializations is essentially Batyrev's $E(X;u,1)$, see [BL2, Proposition 3.7].

\smallskip
(3) We hope that the ideas in this paper can be useful for generalizing also other invariants, for instance the singular elliptic genus, beyond the log terminal case.
\enddemo

\bigskip
\bigskip      
\heading
2. Zeta functions for arbitrary divisors on klt pairs
\endheading
\bigskip
\noindent
{\bf 2.1.}  Let $X$ be a normal variety and $B$ a $\Bbb Q$--divisor on $X$ such that 
$K_X + B$ is $\Bbb Q$--Cartier (when $X$ is a surface we omit this condition).  Let moreover $D$ be any $\Bbb Q$--Cartier divisor on $X$ (again when $X$ is a surface $D$ can be any $\Bbb Q$--divisor).

Take a log resolution $\pi : Y \rightarrow X$ of $(X,\supp B \cup \supp D)$ and denote by $E_i, i \in T,$ the irreducible exceptional divisors of $\pi$ and the strict transforms by $\pi$ of the irreducible components of $\supp B \cup \supp D$.  To each $E_i$ we associate two rational numbers $\nu_i$ and $N_i$, given by
$$\align K_Y = & \pi^\ast (K_X + B) + \sum_{i \in T} (\nu_i - 1)E_i \qquad \text { and }\\
\pi^\ast D = & \sum_{i \in T} N_i E_i ; \endalign$$
hence $\nu_i$ is just the log discrepancy of $E_i$ with respect to the pair $(X,B)$.

\noindent
We use the notations $E_I := \cap_{i \in I} E_i$ and $E^\circ_I := E_I \setminus \cup_{\ell \not \in I} E_\ell$ for $I \subset T$.  In particular $E^\circ_\emptyset = Y \setminus (\cup_{\ell \in T} E_\ell)$, and $Y$ is the disjoint union of the $E^\circ_I, I \subset T$.
\bigskip
\noindent  
{\bf 2.2. Definition.}  When $(X,B)$ is klt we associate to $(X,B)$ and $D$ the `zeta function'
$$\Cal Z(s) = \Cal Z((X,B),D;s) := \sum_{I \subset T} [E^\circ_I] \prod_{i \in I} \frac{L-1}{L^{\nu_i+sN_i}-1} \, .$$
Here $L^{-s}$ should be considered as a variable $T$ and $\frac{L-1}{L^{\nu_i+sN_i}-1}$ as $\frac{(L-1)T^{N_i}}{L^{\nu_i} - T^{N_i}}$ or $\frac{L-1}{L^{\nu_i}T^{|N_i|} -1}$, depending on whether $N_i$ is positive or negative, respectively. So $Z(s)$ lives e.g. in the polynomial ring \lq with fractional powers\rq\ in the variable $T$ over the ring $\Cal R$, localized with respect to the elements $L^{\nu} - T^N$ and $L^\nu T^N -1$ for $\nu \in \Bbb Q_{\geq 0}$ and $N \in \Bbb Q_{>0}$. (Here $\nu>0$ suffices, but we need $\nu=0$ in (2.6).)  

All this may seem weird at first sight, but it is quite similar to the motivic zeta functions and the rings they live in from [DL2] and [Ve2].
\bigskip
\proclaim{2.3. Proposition}  Definition 2.2 does not depend on the chosen resolution.
\endproclaim
\bigskip
\noindent
{\sl Remark.}  We prove this using weak factorization.  See (2.8) for a remark on weak factorization versus motivic integration.  We want to present carefully the main arguments for this independency, since one is frequently sloppy in applying the weak factorization theorem.  In particular part (1') in Theorem 1.7.1 is often crucial, as is the case here.
\bigskip
\noindent
{\bf 2.3.1.} Let $V$ be a smooth irreducible variety and $\cup_{i \in S} F_i$ a normal crossings divisor on $V$ with the $F_i$ irreducible.      
Denote $F^\circ_I := (\cap_{i \in I} F_i) \setminus (\cup_{\ell \not\in I} F_\ell)$ for $I \subset S$.  We associate a zeta function $Z_V(\Cal K, \Cal D; s)$ to two $\Bbb Q$--divisors $\Cal K = \sum_{i \in S} (k_i - 1)F_i$, with all $k_i > 0$, and $\Cal D  = \sum_{i \in S} d_i F_i$ on $V$ :
$$\Cal Z_V (\Cal K, \Cal D;s) := \sum_{I \subset S} [F^\circ_I] \prod_{i \in I} \frac{L-1}{L^{k_i+sd_i}-1},$$
living in a ring as described in (2.2).  In particular,  using the resolution $\pi : Y \rightarrow X$ of (2.1), our proposed definition for $\Cal Z(s)$ in (2.2) is $\Cal Z_Y (K_Y - \pi^\ast (K_X + B), \pi^\ast D; s)$.
\bigskip
\proclaim
{2.3.2. Lemma}  With the notation of (2.3.1), let $h : W \rightarrow V$ be a composition of blowing--ups with smooth centre, having normal crossings with $\cup_{i \in S} F_i$ and its consecutive inverse image.  Then
$$\Cal Z_V (\Cal K, \Cal D; s) = \Cal Z_W(h^\ast \Cal K + (K_W - h^\ast K_V), h^\ast \Cal D; s).$$
\endproclaim
\bigskip
\demo{Proof}  (i) We first suppose that $h$ is just one blowing--up with centre $Z$ of codimension $r \geq 2$ in $V$ and exceptional variety $F$.  Denote the strict transform of $F_i$ in $W$ by  $\tilde F_i$.  Say $Z \subset F_i$ for $1 \leq i \leq m$ (here $0 \leq m \leq r$).  Then 
$$\align
&h^\ast\Cal K + (K_W - h^\ast K_V) = (k-1)F + \sum_{i \in S} (k_i - 1)\tilde F_i \qquad\text{and} \\ 
&h^\ast \Cal D = dF + \sum_{i \in S} d_i \tilde F_i,
\endalign$$
where $k = \sum^m_{i=1} (k_i - 1) + r$ and $d = \sum^m_{i=1} d_i$.  (So indeed $k > 0$, as required in the definition of $\Cal Z_W (h^\ast \Cal K + (K_W - h^\ast K_V), h^\ast \Cal D; s)$.)

We must compare the contribution of $Z$ to $\Cal Z_V(\Cal K, \Cal D; s)$ and the contribution of $F$ to $\Cal Z_W(h^\ast \Cal K + (K_W - h^\ast K_V), h^\ast \Cal D; s)$; i.e. they should be equal.  Since $h |_F : F \rightarrow Z$ is locally a product it is sufficient to compare the contribution of a point $P \in Z$ and of $h^{-1}P(\cong \Bbb P^{r-1}) \subset F$, respectively.  Say $P$ also belongs to $F_j, m  + 1 \leq j \leq n$ (here $m \leq n \leq \dim V$).  Then these contributions are    
$$\prod^n_{i=1} \frac{L-1}{L^{k_i+sd_i}-1} \qquad  \and \qquad \frac{L-1}{L^{k+sd}-1} \cdot \prod^n_{i=m+1} \frac{L-1}{L^{k_i+sd_i}-1} \cdot A^{r-1}_m,$$
respectively, where
$$A^{r-1}_m := \sum_{I \subset \{ 1, \cdots ,m \}} [(h^{-1} P \cap (\cap_{i \in I} \tilde F_i))^\circ] \prod_{i \in I} \frac{L-1}{L^{k_i + sd_i}-1}$$
is the `contribution' of the $m$ hyperplanes $h^{-1} P \cap \tilde F_i$ in general position in $h^{-1}P \cong \Bbb P^{r-1}$.  Here as usual $(h^{-1}P \cap (\cap_{i \in I} \tilde F_i))^\circ$ is $(h^{-1} P \cap (\cap_{i \in I} \tilde F_i)) \setminus \cup_{\ell \in \{ 1 , \cdots , m \} \setminus I} \tilde F_\ell$.

So it is sufficient to show that $\prod^m_{i=1} \frac{L-1}{L^{k_i + sd_i}-1} = \frac{L-1}{L^{k+sd} - 1} \cdot A^{r-1}_m$, or, equivalently, that
$$A^{r-1}_m = (L-1)^{m-1} \frac{L^{\sum^m_{i=1} (k_i - 1) + r + s \sum^m_{i=1} d_i}-1}{\prod^m_{i=1} (L^{k_i + sd_i}-1)}.$$
Using (double) induction on $r$ and $m$, this is easy to verify.
\medskip
(ii) For the general statement it is sufficient to treat the case that $h$ is a composition of two blowing--ups $W_2 \overset h_2 \to \longrightarrow W_1 \overset h_1 \to \longrightarrow V$.  Applying (i) twice we obtain that
$$\split \Cal Z_V (\Cal K, \Cal D; s) & = \Cal Z_{W_1} (h^\ast_1 \Cal K + (K_{W_1} - h^\ast_1 K_V), h^\ast_1 \Cal D; s) \\
& = \Cal Z_{W_2} (h^\ast_2 h^\ast_1 \Cal K + h^\ast_2 K_{W_1} - h^\ast_2 h^\ast_1 K_V + K_{W_2} - h^\ast_2 K_{W_1}, h^\ast_2 h^\ast_1 \Cal D;s) \\
& = \Cal Z_{W_2} (h^\ast \Cal K + K_{W_2} - h^\ast K_V , h^\ast \Cal D; s). \qed
\endsplit$$
\enddemo
\bigskip
\noindent
{\bf 2.3.3.} {\sl Proof of Proposition 2.3.}

\smallskip
\noindent
Let $\pi : Y \rightarrow X$ and $\pi^\prime : Y^\prime \rightarrow X$ be two log resolutions of $(X,\supp B \cup \supp D)$.  By Theorem 1.7.1 the associated birational map $\phi : Y -\! \rightarrow Y^\prime$ (of $X$--varieties) decomposes in a sequence of blowing--ups and blowing--downs as in Figure 1.     

Here $\phi = h_{2m} \circ h^{-1}_{2m-1} \circ h_{2m-2} \circ \cdots \circ h^{-1}_3 \circ h_2 \circ h^{-1}_1$, and all $h_{2i-1} : Y_{2i-1} \rightarrow Y_{2i-2}$ and $h_{2i} : Y_{2i-1} \rightarrow Y_{2i}$ are compositions of blowing--ups with smooth centres. (Some  $h_j$ can be the identity, making the argument just easier.) 
The $\lambda_{2i}$ and $\mu_{2i}$ are morphisms, see (1.7.1(1')), and 
furthermore $\alpha := \pi \circ \lambda_{2k} = \pi' \circ \mu_{2k}$.  

By (1.7.1(2)) we may suppose that we have the following normal crossing properties.  For $i \leq k$ the inverse images under $\lambda_{2i} : Y_{2i} \rightarrow Y$ of the union of the exceptional divisor of $\pi$ and the strict transform in $Y$ of $\supp B \cup \supp D$ are also normal crossings divisors in $Y_{2i}$, and the centres of the blowing--ups occurring in the $h_{2i-1}$ and $h_{2i}$ have normal crossings with their consecutive inverse images.  For $i \geq k$ we have the analogous statements for the $\mu_{2i}$, and for the $h_{2i+1}$ and $h_{2i+2}$.

\vskip 1truecm
\centerline{\beginpicture
\setcoordinatesystem units <.5truecm,.5truecm>
\put {$Y_1$} at  -9.5 4.7
\put {$Y_3$} at -6 4.7
\put {$Y_{2k-1}$} at -1.5 4.7
\put {$Y_{2k+1}$} at 1.5 4.7
\put {$Y_{2m-3}$} at 6 4.7
\put {$Y_{2m-1}$} at 9.5 4.7
\arrow <.3truecm> [.2,.6] from 1.5 4 to .5 1
\arrow <.3truecm> [.2,.6]  from 9 4 to 8 1
\arrow <.3truecm> [.2,.6]  from -6 4 to -7 1
\arrow <.3truecm> [.2,.6] from -1.5 4 to -.5 1
\arrow <.3truecm> [.2,.6]  from -9 4 to -8 1
\arrow <.3truecm> [.2,.6]  from 6 4 to 7 1
\arrow <.3truecm> [.2,.6] from 10 4 to 13 -1
\arrow <.3truecm> [.2,.6]  from -10 4 to -13 -1
\put {$h_1$} at  -12 2
\put {$h_2$} at  -9 2.5
\put {$h_3$} at  -5.9 2.5
\put {$h_{2k}$} at  -1.6 2.5
\put {$h_{2k+1}$} at  1.3 2.5
\put {$h_{2m-2}$} at  6 2.5
\put {$h_{2m-1}$} at  9 2.5
\put {$h_{2m}$} at  12.2 2
\put {$Y_2$} at  -7.5 .2
\put {$Y_{2k}$} at 0 -.2
\put {$Y_{2m-2}$} at 7.5 .2
\arrow <.3truecm> [.2,.6] from -8.5 -.5 to -12 -2
\arrow <.3truecm> [.2,.6]  from -1.5 -.5 to -12 -3
\arrow <.3truecm> [.2,.6] from 8.5 -.5 to 12 -2
\arrow <.3truecm> [.2,.6]  from 1.5 -.5 to 12 -3
\put {$\lambda_2$} at  -10.5 -.6
\put {$\lambda_{2k}$} at  -6 -2.3
\put {$\mu_{2k}$} at  6 -2.1
\put {$\mu_{2m-2}$} at  10.4 -.7
\put {$Y=Y_0$} at  -14 -2.2
\put {$X$} at 0 -6
\put {$Y_{2m}=Y'$} at 14 -2.2
\arrow <.3truecm> [.2,.6]  from -13 -3.5 to -1.5 -6
\arrow <.3truecm> [.2,.6] from 13 -3.5 to 1.5 -6
\arrow <.3truecm> [.2,.6]  from 0 -1.5 to 0 -5
\put {$\alpha$} at  .5 -3
\put {$\pi$} at  -6.5 -5.5
\put {$\pi'$} at  6.5 -5.5
\multiput {$\dots$} at -4 2.5  4 2.5  -9 -1.5  9 -1.5 /
\endpicture}
\vskip 1truecm
\centerline{\smc Figure 1}
\vskip 1truecm

With the concept of the zeta function from (2.3.1) we have to show that
$$\Cal Z_Y (K_Y - \pi^\ast (K_X + B), \pi^\ast D;s) = \Cal Z_{Y^\prime} (K_{Y^\prime} - \pi^{\prime \ast} (K_X + B), \pi^{\prime \ast} D;s).$$
We abbreviate $\Cal K := K_Y - \pi^\ast (K_X + B)$ and $\Cal D := \pi^\ast D$.  Then 
$$\split
\Cal Z_Y(\Cal K, \Cal D;s) & = \Cal Z_{Y_1} (h^\ast_1 \Cal K + K_{Y_1} - h^\ast_1 K_Y, h^\ast_1 \Cal D; s) \\
& = \Cal Z_{Y_1} (h^\ast_2 \lambda^\ast_2 \Cal K + K_{Y_1} - h^\ast_2 K_{Y_2} + h^\ast_2 (K_{Y_2} - \lambda^\ast_2 K_Y), h^\ast_2 \lambda^\ast_2 \Cal D; s) \\ & = \Cal Z_{Y_2} (\lambda^\ast_2 \Cal K + K_{Y_2} - \lambda^\ast_2 K_Y, \lambda^\ast_2 \Cal D; s),
\endsplit
$$    
where the first and third identities are Lemma 2.3.2, applied to $h_1$ and $h_2$, respectively, and the second one is straightforward, but requires that $\lambda_2$ is a morphism !  Proceeding further analogously we obtain that
$$\split
\Cal Z_Y (\Cal K, \Cal D; s) & = \Cal Z_{Y_2} (\lambda^\ast_2 \Cal K + K_{Y_2} - \lambda^\ast_2 K_Y, \lambda^\ast_2 \Cal D; s) \\
& = \Cal Z_{Y_4} (\lambda^\ast_4 \Cal K + K_{Y_4} - \lambda^\ast_4 K_Y, \lambda^\ast_4 \Cal D; s) \\
& = \cdots = \Cal Z_{Y_{2k}} (\lambda^\ast_{2k} \Cal K + K_{Y_{2k}} - \lambda^\ast_{2k} K_Y, \lambda^\ast_{2k} \Cal D; s).
\endsplit
$$ 
By definition of $\Cal K$ and $\Cal D$ this is just
$$\split \Cal Z_{Y_{2k}} & (\lambda^\ast_{2k} K_Y - \lambda^\ast_{2k} \pi^\ast (K_X + B) + K_{Y_{2k}} - \lambda^\ast_{2k} K_Y, \lambda^\ast_{2k} \pi^\ast D; s) \\ & = \Cal Z_{Y_{2k}} (K_{Y_{2k}} - \alpha^\ast (K_X + B), \alpha^\ast D; s). \endsplit$$ 
Now completely analogously we see that also
$$\Cal Z_{Y^\prime}  (K_{Y^\prime} - \pi^{\prime \ast} (K_X + B), \pi^{\prime \ast} D;s) =  \Cal Z_{Y_{2k}} (K_{Y_{2k}} - \alpha^* (K_X + B), \alpha^\ast D; s). \qed
$$
\bigskip
\noindent
{\bf 2.4. Definition.} When $(X,B)$ is klt we associate to $(X,B)$ and $D$ the zeta function
$$z(s) = z((X,B),D;s) := \sum_{I \subset T} \chi (E^\circ_I) \prod_{i \in I} \frac{1}{\nu_i + sN_i} \in \Bbb Q(s).$$
\bigskip
\noindent
{\bf 2.4.1.} {\sl Remarks.}
 (i) This zeta function can be seen as a specialization of $\Cal Z(s)$; see [DL2] or [Ve2].

(ii) For $z(s)$ the independency of the chosen resolution can also be derived in an  elementary way from the $e$--invariant for effective divisors on $X$ [Ve2].  In fact we need more generally such an $e$--invariant associated to an effective divisor and the pair $(X,B)$, but this generalization is straightforward.

Set $D = D^+ - D^-$, where $D^+$ and $D^-$ are effective and have no common component, and decompose accordingly $N_i = N^+_i - N^-_i$ for $i \in T$.  (To be precise we must restrict here to the case that $D^+$ and $D^-$ are $\Bbb Q$--Cartier; but anyhow this is a preparation for \S 3 where we will use the zeta functions of this section for a $\Bbb Q$--factorial $X$.)  Consider now $e(m_1 D^+ + m_2 D^-)$ for all $m_1,m_2 \in  \Bbb Z_{\geq 0}$.  The function
$$z^\prime (D^+,D^-; s_1,s_2) := \sum_{I \subset T} \chi(E^\circ_I) \prod_{i \in I} \frac{1}{\nu_i + s_1 N^+_i + s_2 N^-_i} \in \Bbb Q (s_1,s_2)$$
is the unique rational function in two variables $s_1,s_2$ yielding $e(m_1D^+ + m_2 D^-)$ when evaluating $s_1$ in $m_1$ and $s_2$ in $m_2$.

Then $z(s) = z^\prime (D^+,D^-;s,-s)$.
\bigskip
\noindent
{\bf 2.5.} (i) One can of course introduce analogously `intermediate level' zeta functions to a klt pair $(X,B)$ and a divisor $D$ on $X$, e.g. on the level of Hodge polynomials.  As long as the coefficient ring has no zero divisors, the argument in 2.4.1(ii) should work.

(ii) For any constructible subset $W$ of $X$, we can introduce more generally zeta functions $\Cal Z_W(s) = \Cal Z_W ((X,B);s)$ using $[E^\circ_I \cap \pi^{-1} W]$ instead of $[E^\circ_I]$ in Definition 2.2.  (This is also true for the invariants in [DL2] and [Ve2].)  Some interesting cases are $W = X_{\sing}, W = B$, and in particular $W = \{ P \}$ for some point $P \in X$.  Then we rather write $\Cal Z_P(s)$; this is the appropriate invariant when studying singularity germs.

Of course one can treat the zeta functions on other levels (and the next ones in (2.6)) also in this more general $W$--setting.
\bigskip
\noindent
{\bf 2.6.} From a pragmatic point of view one can define $\Cal Z(s)$ and its specializations for klt pairs and arbitrary divisors because all $\nu_i \ne 0$ in any log resolution.  In fact one can define such zeta functions as long as in suitable log resolutions either $\nu_i \ne 0$ or $N_i \ne 0$.  We will need the following case.
\bigskip
\noindent
{\bf 2.6.1. Definition.}  Let $(X,B = \sum_i b_i B_i)$ be a dlt pair and let $Z \subset X$ be the closed subset of Definition 1.4.1 (iii), where $B|_{X \setminus Z}$ is a reduced normal crossings divisor.  Here we only consider $\Bbb Q$--divisors $D$ on $X$ such that $\supp D|_{X \setminus Z} = \supp B|_{X \setminus Z}$, and the coefficients of $D|_{X \setminus Z}$ are either all positive or all negative.

Then we define $\Cal Z(s)$ and $z(s)$ as in (2.2) and (2.4).
\bigskip
\noindent
{\bf 2.6.2.}  This definition is independent of the chosen  resolution by the same arguments as before; we just have to check for any log resolution as in (2.1) that for $i \in T$ either $\nu_i \ne 0$ or $N_i \ne 0$.

If $E_i$ is exceptional and $\pi(E_i) \subset Z$, then $\nu_i > 0$.  Also, if $E_i$ is the strict transform of a $B_i$ with $b_i < 1$, then $\nu_i = 1 - b_i > 0$.  If $E_i$ is another exceptional component, or the strict transform of a $B_i$ with $b_i = 1$, then $\nu_i \geq 0$, and when $\nu_i = 0$, then $N_i \ne 0$ by our assumption on $D$.
\bigskip
\noindent
{\bf 2.7.} The stringy $E$--function $E(X)$ of a projective log terminal variety $X$ satisfies a `Poincar\'e duality' result [Ba1, Theorem 3.7] :
$$(uv)^{\dim X} E(X)|_{\Sb u \rightarrow u^{-1} \\ v \rightarrow v^{-1} \endSb} = E(X).$$      
This generalizes to a `functional equation' for the zeta functions introduced in (2.2) and (2.6.1), specialized to the level of Hodge polynomials.  With the notation of (2.1) these zeta functions are of the form
$$Z((X,B),D;s) := \sum_{I \subset T} H(E^\circ_I) \prod_{i \in I} \frac{uv - 1}{(uv)^{\nu_i + sN_i}-1},$$
where now $(uv)^{-s}$ should be considered as a variable $T$.
\bigskip
\proclaim{2.7.1. Proposition}  Let $(X,B)$ and $D$ be as in either (2.2) or (2.6.1), and let moreover $X$ be projective.  Then
$$(uv)^{\dim X} Z((X,B),D;s)|_{\Sb u \rightarrow u^{-1} \\ v \rightarrow v^{-1} \endSb} = Z((X,B),D;s).$$
This substitution has to be interpreted literally : the variable $T = (uv)^{-s}$ must be replaced by $T^{-1} = (uv)^s$; so in particular $(uv)^{\nu_i + sN_i}$ is replaced by $(uv)^{-\nu_i - sN_i}$.
\endproclaim
\bigskip
\demo{Proof}  Analogously as in [Ba1, Theorem 3.7] or in [DM] one easily sees that an alternative expression for $Z((X,B),D;s)$ is
$$\sum_{I \subset T} H(E_I) \prod_{i \in I} \Big(\frac{uv-1}{(uv)^{\nu_i+sN_i}-1} - 1 \Big).$$
Then just as in [Ba1] and [DM] the stated functional equation is true because it is valid for each term in the above sum, using ordinary Poincar\'e duality for the smooth projective varieties $E_I$. 
\qed \enddemo

\bigskip
\noindent
{\bf 2.7.2.} {\sl Remark.} Using the duality involution from [Bi, Corollary 3.4], Batyrev's Poincar\'e duality and more generally our Proposition 2.7.1 can be \lq upgraded\rq\ to the level of the Grothendieck ring. (Then  we have to redefine the ring $\Cal R$ to make it contain $L^{-1}$.)

\bigskip
\noindent
{\bf 2.8. Remark.}  On the level of the Grothendieck ring one obtains a priori finer invariants using weak factorization than with motivic integration.  For example when we take $D=0$ in (2.2), we just obtain for the klt pair $(X,B)$ a Batyrev--type invariant
$$\Cal E(X,B) = \sum_{I \subset T} [E^\circ_I] \prod_{i \in I} \frac{L-1}{L^{a_i}-1},$$
using the notation of (0.1), which specializes to the stringy $E$--function $E(X,B)$.  It lives in the ring $\Cal R$ of (1.8(i)); this $\Cal R$ was a subring of a localization of $K_0(\Var_{\Bbb C})$, extended with fractional powers of $L$.  Alternatively, one can introduce such an invariant, using motivic integration techniques as in [Ba2] and [DL3].  But then we only know that it lives in some completion of $K_0(\Var_{\Bbb C})$ (extended with fractional powers of $L$), more precisely in the image of $\Cal R$ in this completion.  Since it is not known whether the map of $\Cal R$ to this completion is injective, the first invariant $\Cal E(X,B) \in \Cal R$ is a priori finer.  

We illustrate this in the case that $K_X + B$ is Cartier.  Then all log discrepancies $a_i$ are positive integers, and we can consider $\Cal E(X,B)$ simply in the localization $\Cal R^\prime$ of $K_0(\Var_{\Bbb C})$ with respect to $[\Bbb P^{a-1}] = L^{a-1} + \cdots + L + 1, a \in \Bbb Z_{> 1}$.  With motivic integration one considers the completion of $K_0(\Var_{\Bbb C})[L^{-1}]$ with respect to the decreasing filtration by subgroups $F^m$, generated by the elements $[S]L^{-i}$ with $\dim S- i \leq -m$.  And then the analogous invariant lives in the image of $\Cal R^\prime$ in this completion. 

See also [DL3] and [Ve2] for a description of such completions.  Analogously, one obtains for some invariants in [DL3] and [Ve2], where  motivic integration is used, a priori finer ones using weak factorization.

\bigskip
\noindent
{\bf 2.8.1.} In fact, Kontsevich introduced motivic integration to prove that birationally equivalent (smooth, complete) Calabi--Yau varieties have the same Hodge numbers. More precisely, he showed that these varieties induce the same element in the above described completion of $K_0(\Var_{\Bbb C})[L^{-1}]$, see [DL4, 4.4.2]. By essentially the same arguments as for Proposition 2.3, one shows the following finer result.

\bigskip
\proclaim{Proposition} Let $Y$ and $Y'$ be birationally equivalent complete smooth Calabi--Yau varieties (i.e. $K_Y=K_{Y'}=0$). Then $[Y]=[Y']$ in the the ring $\Cal R'$.
\endproclaim

\bigskip
\bigskip
\head
3. Stringy zeta functions, generic case
\endhead
\bigskip
\noindent
{\bf 3.1.} Let $X$ be a normal variety and $B = \sum_i b_i B_i$ a $\Bbb Q$--divisor on $X$, where the $B_i$ are distinct and irreducible and all $b_i$ satisfy $0 \leq b_i < 1$, such that $K_X + B$ is $\Bbb Q$--Cartier.  (When $X$ is a surface we omit this last condition.)  

Take a (relative) log minimal model $p : X^m \rightarrow X$ of $(X,B)$ as defined in (1.6.1(1)); we denote by $F = \sum_{i \in T^m} F_i$ the (reduced) exceptional divisor of $p$, and by $B^m$ the strict transform of $B$ by $p$.  We assume to be in the generic case that all log discrepancies with respect to $(X,B)$ of exceptional divisors of $p$ are negative, i.e.
$$K_{X^m} + B^m + F = p^\ast (K_X + B) + \sum_{i \in T^m} a_i F_i$$
with all $a_i < 0$.
In (3.7) we will explain that this condition is indeed generic and conceptual : it is equivalent to just asking that $(X,B)$ has no strictly lc singularities.
\bigskip
\noindent
{\bf 3.2. Definition.}  We assume the (relative log) MMP.  To $(X,B)$ as in (3.1) we associate the {\it stringy zeta function}
$$\Cal Z(s) = \Cal Z((X,B);s) := \Cal Z((X^m,B^m + F), (K_{X^m} + B^m  + F) - p^\ast (K_X + B);s), $$
where the right hand side is the zeta function associated in (2.6.1) to the dlt pair $(X^m,B^m + F)$ and the \lq log discrepancy divisor\rq\ $D :=\sum_{i \in T^m} a_i F_i = K_{X^m} + B^m + F  - p^\ast (K_X + B)$ on $X^m$. 
(When $B=0$ we just write $\Cal Z(X;s)$.)

We check that the condition in (2.6.1) on $\supp D$ is satisfied.  Let $Z \subset X^m$ denote the closed subset of (1.4.1(iii)), where $(B^m + F)|_{X^m \setminus Z}$ is a reduced normal crossings divisor. Note that then $(B^m + F)|_{X^m \setminus Z}= F|_{X^m \setminus Z}$ because all $b_i<1$. Since $D = \sum_{i \in T^m} a_i F_i$ where by assumption all $a_i < 0$, we have indeed that $\supp D|_{X^m \setminus Z} = \supp F|_{X^m \setminus Z}= \supp(B^m + F)|_{X^m \setminus Z}$, and the coefficients of $D|_{X^m \setminus Z}$ are all negative.
\bigskip
\noindent
{\bf 3.2.1. Formula.}  Take a log resolution $h : Y \rightarrow X^m$ of the pair $(X^m,B^m + F)$.  Let $E_i, i \in T$, be the irreducible components of the exceptional divisor of $h$ and of the strict transform of $B^m + F$.  In particular $E_i, i \in T^m \subset T$, is the strict transform of $F_i$ in $Y$.  As usual we denote $E^\circ_I := (\cap_{i \in I} E_i) \setminus (\cup_{\ell \notin I} E_\ell)$ for $I \subset T$.  Then
$$\Cal Z(s) = \sum_{I \subset T} [E^\circ_I] \prod_{i \in I} \frac{L - 1}{L^{\nu_i  + sN_i} - 1} ,$$
where
$$\align K_Y &= h^\ast (K_{X^m} + B^m + F) + \sum_{i \in T} (\nu_i - 1)E_i
\qquad\qquad \text{and} \\
h^\ast (D) &= h^\ast (K_{X^m} + B^m + F - p^\ast (K_X + B)) = \sum_{i \in T} N_i E_i.
\endalign$$
\bigskip
\proclaim
{3.2.2. Proposition}  Definition 3.2 does not depend on the chosen log minimal model.
\endproclaim
\bigskip
\demo{Proof}  Let $p_1 : X_1 \rightarrow X$ and $p_2 : X_2 \rightarrow X$ be two log minimal models of $(X,B)$ as in (3.1).  We denote for $j = 1,2$ by $F^j$ and $B^j$ the exceptional divisor of $p_j$ and the strict transform of $B$ in $X_j$, respectively.

Take a common log resolution $Y$ of the $(X_j,B^j + F^j)$ as in Figure 2, and denote by $\tilde B$ the strict transform of $B$ in 
$Y$.  We have that $h^\ast_1 (K_{X_1} + B^1 + F^1) = h^\ast_2 (K_{X_2} + B^2 + F^2)$; see the proof of [KM, Theorem 3.52(2)].  Hence 
$$\align &K_Y - h^\ast_1 (K_{X_1} + B^1 + F^1) = K_Y - h^\ast_2 (K_{X_2} + B^2 + F^2) \qquad \text{and} \\
&h^\ast_1 (K_{X_1}  + B^1 + F^1 - p^\ast_1 (K_X + B)) = h^\ast_2 (K_{X_2} + B^2 + F^2 - p^\ast_2 (K_X + B)).
\endalign$$
Using Formula 3.2.1, this means that indeed $X_1$ and $X_2$ yield the same right hand side in (3.2). \qed
\enddemo
\vskip 1true cm
\centerline{\beginpicture
\setcoordinatesystem units <.5truecm,.5truecm>
\put{$\searrow$} at 1 1
\put{$\swarrow$} at 1 -1
\put{$\searrow$} at -1 -1
\put{$\swarrow$} at -1 1
\put{$h_1$} at -1.5 1.5
\put{$h_2$} at 1.5 1.5
\put{$p_2$} at 1.5 -1.5
\put{$p_1$} at -1.4 -1.4
\put{$X_1$} at -2 0
\put{$Y$} at  0 2
\put{$X_2$} at 2 0
\put{$X$} at  0 -2
\endpicture}
\vskip 1truecm
\centerline{\smc Figure 2}
\vskip 1truecm
\noindent
{\bf 3.2.3} (i) To $(X,B)$ as in (3.1) we can associate analogously stringy zeta functions on other levels, e.g. on the level of Euler characteristics.

(ii) For a constructible subset $W$ of $X$ we can introduce more generally $\Cal Z_W(s) = \Cal Z_W((X,B);s)$ as
$$\Cal Z_{p^{-1}W} ((X^m,B^m + F), (K_{X^m} + B^m + F) - p^{\ast} (K_X + B);s),$$
see (2.5.(ii)).  Then in formula (3.2.1) one must replace 
$[E^\circ_I]$ by $[E^\circ_I \cap (p\circ h)^{-1} W]$.  The same remark applies to other levels.
\bigskip
\noindent
{\bf 3.3.} When the pair $(X,B)$ is itself klt and $p : X^m \rightarrow X$ is a log minimal model of $(X,B)$, then $p$ has no exceptional divisors (see e.g. (3.6.1)).  In particular the generic condition in (3.1) is trivially verified.  So $\Cal Z((X,B);s)$ is the zeta function of (2.2) associated to $(X^m,B^m)$ and the divisor $D = 0$.  Choosing $h$ in (3.2.1) such that $\pi = p \circ h$ is a log resolution of $(X,B)$ we see that $\Cal Z(s)$ is just the stringy $\Cal E$--invariant $\Cal E(X,B)$ of Batyrev. 
\bigskip
\noindent
{\bf 3.4.}  We now make the link with question (I) in the introduction.  We consider a pair $(X,B)$ as in (3.1) which has some log resolution $\pi : Y \rightarrow X$ for which {\it all} log discrepancies with respect to $(X,B)$ of divisors on $Y$ are nonzero; i.e.
$$K_Y = \pi^\ast(K_X + B) + \sum_{i \in T} (a_i - 1)E_i \quad \text { with all } a_i \ne 0,$$    
where as usual $E_i, i \in T$, are the irreducible exceptional divisors of $\pi$ and the strict transforms by $\pi$ of the irreducible components of $\supp B$.

Assume now that $\pi$ factors through a log minimal model $p : X^m \rightarrow X$ of $(X,B)$, i.e.
$$\pi : Y  \overset h \to \longrightarrow X^m \overset p \to \longrightarrow X, $$
where $h$ is a {\it morphism}. 
\bigskip
\noindent
{\bf 3.4.1. Claim.}  $\sum_{I \subset T} [E^\circ_I] \prod_{i \in I} \frac{L-1}{L^{a_i}-1} = \Cal Z((X,B);1)$.
\bigskip
\demo{Remark} Be evaluating $\Cal Z((X,B);s)$ in $1$ we mean evaluating $T=L^{-s}$ in $L^{-1}$. One verifies that this yields a well defined element in the ring $\Cal R$ of (1.8(i)).
\enddemo
\bigskip
\demo{Proof}  By definition $\Cal Z((X,B);s)$ is determined by the $\nu_i$ and $N_i$ in
$$\align &K_Y = h^\ast (K_{X^m} + B^m + F) + \sum_{i \in T} (\nu_i - 1)E_i
\qquad\qquad \text{and} \\
&h^\ast (K_{X^m} + B^m + F - p^\ast (K_X + B)) = \sum_{i \in T} N_i E_i, 
\endalign$$
using the notation of (3.1).  Adding both equalities yields $K_Y = \pi^\ast (K_X + B) + \sum_{i \in T} (\nu_i + N_i - 1)E_i$ and thus $\nu_i + N_i = a_i \ne 0$ for all $i \in T$.  So the evaluation $\Cal Z((X,B);1)$ indeed makes sense and is as stated. \qed
\enddemo
\bigskip
\noindent
{\bf 3.4.2.}  We conclude that the expressions
$$\sum_{I \subset T} [E^\circ_I] \prod_{i \in I} \frac{L-1}{L^{a_i}-1}$$
are the same for all log resolutions (with all $a_i \ne 0$) {\it that factorize through some log minimal model} of $(X,B)$; here this model can depend on the resolution.  In this restricted sense they can be considered as a generalized $\Cal E$--invariant.
\bigskip
\noindent
{\bf 3.5.}  In [Ve3] we associated stringy invariants to {\it any} normal surface $X$ without strictly lc singularities.  We recall their definition, but first we state the structure theorem  on which it is based.
\bigskip
\noindent
\proclaim{3.5.1. Theorem} {\rm [Ve3, 2.10]}  Let $P \in X$ be a normal surface singularity germ which is not log canonical.  Let $\pi : Y \rightarrow X$ be the minimal log resolution of $P \in X$; denote the irreducible components of $\pi^{-1}P$ by $E_i, i \in T$, and their log discrepancies with respect to $X$ by $a_i$.  Then $\pi^{-1} P = \cup_{i \in T} E_i$ consists of the connected part $\Cal N = \cup_{i \in T, a_i <0} E_i$, to which a finite number of chains are attached as in Figure 3.  Here $E_i \subset \Cal N, E_\ell \cong \Bbb P^1$ for $1 \leq \ell \leq r, a_1 \geq 0$ and $(a_i <) a_1 < a_2 < \cdots < a_r < 1$.  

In particular, if $a_j=0$, then $E_j \cong \Bbb P^1$ and $E_j$ intersects exactly one or two other components $E_i$ (and those have $a_i \ne 0$).
\endproclaim

\vskip 1true cm
\centerline{\beginpicture
\setcoordinatesystem units <.5truecm,.5truecm>
\putrule from -4 3 to 1.5 3
\setlinear  \plot  -1.33 2  4 6 /          \plot 2 6  5 2    /
            \plot  3 2  5.66 4 /         \plot 11 6  14 2 /
            \plot  10.33 4  13 6 /
            \plot  12 2  17.33 6 /        
 \setdashes  \plot   7 5  5.66 4  /       
             \plot   9 3  10.33 4 /      
\putrule from -6 3 to -4 3
 \setsolid
\put {\dots} at 8 4
\put {$E_r$} at 16.9 4.9
\put {$E_{r-1}$} at 13.4 4.5
\put {$E_{2}$} at 4 4.5
\put {$E_1$} at .9 4.5
\put {$E_i$} at -3 3.6
\endpicture}
\vskip 1truecm
\centerline{\smc Figure 3}
\vskip 1truecm
\noindent
{\bf 3.5.2.}  As a corollary, the condition that a normal surface $X$ has no strictly lc singularities is equivalent to the condition in (3.1) that all log discrepancies with respect to $X$ of exceptional divisors on a log minimal model of $X$ are negative.  Indeed, this can be derived from the well known fact that chains $\cup^r_{\ell = 1} E_\ell$ as above are contracted while constructing a log minimal model, starting from $Y$.
\bigskip
\noindent
{\bf 3.5.3. Definition} [Ve3, 3.4].  Let $X$ be a normal surface without strictly lc singularities.  Let $\pi : Y \rightarrow X$ be the minimal log resolution of $X$ and $E_i, i \in T$, the irreducible exceptional curves of $\pi$ with log discrepancy $a_i$ with respect to $X$.  As usual we put $E^\circ_I := (\cup_{i \in I} E_i) \setminus (\cup_{\ell \not \in I} E_\ell)$ for $I \subset T$.  Here we denote also $Z := \{ i \in T \mid a_i = 0 \}$ and $-\kappa_i = E^2_i$ for $i \in Z$.

The stringy $\Cal E$--invariant and stringy Euler number of $X$ are
$$\Cal E (X) := \sum_{I \subset T \setminus Z} [E^\circ_I] \prod_{i \in I} \frac{L-1}{L^{a_i}-1} + \sum_{i \in Z} \frac{\kappa_i(L-1)^2}{(L^{a_{i_1}}-1)(L^{a_{i_2}}-1)}$$
and
$$e(X) := \sum_{I \subset T \setminus Z} \chi(E^\circ_I) \prod_{i \in I} \frac{1}{a_i} + \sum_{i \in Z} \frac{\kappa_i}{a_{i_1}a_{i_2}},$$
respectively, where $E_i, i \in Z$, intersects either $E_{i_1}$ and $E_{i_2}$ or only $E_{i_1}$ (and then we put $a_{i_2} := 1$).
Remark that, with the notation of (1.8(iii)), $\chi(\Cal E(X)) = e(X)$.
\bigskip
This way of dealing with zero log discrepancies is in fact quite natural, see [Ve3] for motivation and results.
Here we also want to mention a recent result of N\'emethi and Nicolaescu in the last part of [NN], where they study weighted homogeneous (hyper)surface singularities. In some Taylor expansion associated to those singularities this generalized $e(X)$ appears, yielding a topological interpretation of it. 

We now show that for normal surfaces $X$ the stringy zeta function $\Cal Z(X;s)$ of (3.2) specializes to $\Cal E(X)$, confirming the naturality of both definitions.
\bigskip
\proclaim{3.5.4. Proposition}  Let $X$ be a normal surface without strictly lc singularities.  Then $\Cal E (X) = \Cal Z(X;1)$.
\endproclaim
\bigskip
\noindent
{\sl Remark.}  A priori it is not clear that the evaluation $\Cal Z(X;1)$ makes sense!  Compare with (3.4.1); now some $a_i$ really can be zero.
\bigskip
\noindent
\demo{Proof}  We use the notation of (3.5.3).  Consider the factorization
$$\pi : Y \overset h \to \longrightarrow X^m \overset p \to \longrightarrow X$$
of $\pi$ through a log minimal model $p : X^m \rightarrow X$ of $X$.  As usual we denote by $F$ the reduced exceptional divisor of $p$,
\CenteredTagsOnSplits
$$\split
& K_Y   = h^\ast(K_{X^m} + F) + \sum_{i \in T} (\nu_i - 1)E_i \qquad \text{and} \\
& h^\ast (K_{X^m} + F -  p^\ast K_X) = \sum_{i \in T} N_i E_i.
\endsplit
 \tag $\dagger$  $$
\TopOrBottomTagsOnSplits
Recall from the proof of (3.4.1) that $\nu_i + N_i = a_i$ for $i \in T$.  

Each $I \subset T \setminus Z$ contributes a term
$$[E^\circ_I] \prod_{i \in I} \frac{L-1}{L^{a_i}-1} \qquad \text{and} \qquad [E^\circ_I] \prod_{i \in I} \frac{L-1}{L^{\nu_i+sN_i}-1}$$
to $\Cal E(X)$ and $\Cal Z(X;s)$, respectively.  So clearly its contributions to $\Cal E (X)$ and $\Cal Z(X;1)$ are the same.  It remains to verify equality of the remaining terms.

Fix $i \in Z$; thus $a_i = 0$ and $a_{i_1} \ne 0 \ne a_{i_2}$.  We must show that the evaluation in $s=1$ of
$$\frac{L-1}{L^{\nu_i+sN_i}-1} \bigl(L-1 + \frac{L-1}{L^{\nu_{i_1} + sN_{i_1}} - 1} + \frac{L-1} {L^{\nu_{i_2} + sN_{i_2}} - 1} \bigr) \tag $*$ $$
equals
$$\frac{\kappa_i(L-1)^2}{(L^{a_{i_1}}-1)(L^{a_{i_2}}-1)}\, . \tag $**$ $$     
It is important here that $E_i$ is exceptional for $h : Y \rightarrow X^m$.  (Maybe $E_{i_1}$  is not exceptional.)  Then intersecting with $E_i$ in $(\dagger)$ and the adjunction formula yield $\kappa_i \nu_i = \nu_{i_1} + \nu_{i_2}$ and $\kappa_i N_i = N_{i_1} + N_{i_2}$.  This implies that $(\ast)$ is equal to
$$\align &  \frac{(L-1)^2}{L^{\nu_i+sN_i}-1}  \cdot \frac{L^{\nu_{i_1} + \nu_{i_2}+s(N_{i_1} + N_{i_2})} - 1}{(L^{\nu_{i_1} + sN_{i_1}} - 1)(L^{\nu_{i_2} + sN_{i_2}} - 1)}\\
& = \frac{(L-1)^2 (L^{\kappa_i(\nu_i + sN_i)} - 1)}{(L^{\nu_i + sN_i} - 1)(L^{\nu_{i_1} + sN_{i_1}} - 1)(L^{\nu_{i_2} + sN_{i_2}} - 1)} \\
& = \frac{(L-1)^2 \sum^{\kappa_i-1}_{j=0} L^{j(\nu_i + sN_i)}} {(L^{\nu_{i_1} + sN_{i_1}} - 1)(L^{\nu_{i_2} + sN_{i_2}} - 1)} \, .
\endalign
$$ 
So, indeed, evaluating ($\ast$) in $s=1$ yields ($\ast \ast$). \qed \enddemo 

\bigskip
\noindent
{\bf 3.6. Example.} Here we mention a concrete example of a threefold singularity $P \in X$, having an exceptional surface with log discrepancy zero in a log resolution, and such that nevertheless 
$\lim_{s\to 1} z_P(X;s) \in \Bbb Q$, i.e. such that the evaluation
$z_P(X;1)$ makes sense.

\smallskip
Let $X$ be the hypersurface $\{x^4+y^4+z^4+t^5=0\}$ in $\Bbb A^4$; its only singular point is $P=(0,0,0,0)$. We sketch the following constructions in Figure 4; we denote varieties and their strict transforms by the same symbol. 

The blowing--up $\pi_1:Y_1\to X$ with centre $P$ is already a resolution of $X$ ($Y_1$ is smooth). Its exceptional surface $E_1$ is the affine cone over the smooth projective plane curve $C=\{x^4+y^4+z^4=0\}$. Let $\pi_2:Y_2\to Y_1$
be the blowing--up with centre the vertex $Q$ of this cone, and exceptional surface $E_2 \cong \Bbb P^2$. Then $E_1 \subset Y_2$ is a ruled surface over $C$ which intersects $E_2$ in a curve isomorphic to $C$. The composition $\pi=\pi_1 \circ \pi_2$ is a log resolution of $P\in X$, and one easily verifies that the log discrepancies are $a_1=0$ and $a_2=-1$; in particular $P\in X$ is not log canonical.

Now $E_1 \subset Y_2$ can be contracted (more precisely one can check that the numerical equivalence class of the fibre of the ruled surface $E_1$ is an extremal ray). Let $h:Y_2\to X_m$ denote this contraction, and let $\pi=p\circ h$. As the notation suggests, one can verify that $K_{X^m}+E_2$ is $p$--nef, implying that $(X^m, E_2)$
is a relative log minimal model of $P\in X$. 

\vskip 1true cm
\centerline{\beginpicture
\setcoordinatesystem units <.45truecm,.45truecm> 
\putrectangle corners at 6 14 and 16 4
\multiput {$E_2$} at 15 7.7  25 2 /
\multiput {$E_1$} at 8.3 13   1.6 4.2  /
\multiput {$C$} at 11.5 8.8   23.5 -.2  /
\ellipticalarc axes ratio 3:1.1  360 degrees from 13 12 center at 11 12
\ellipticalarc axes ratio 3:1.1  360 degrees from 13 9 center at 11 9
\ellipticalarc axes ratio 3:1.1  360 degrees from 15 9 center at 11 9
\ellipticalarc axes ratio 3:1.1  360 degrees from 13 6 center at 11 6
\putrule from 13 6 to 13 12
\putrule from 9 6 to 9 12
\putrectangle corners at -3 -5 and 3 5
\ellipticalarc axes ratio 3:1.1  360 degrees from 2 3 center at 0 3
\ellipticalarc axes ratio 3:1.1  360 degrees from 2 -3 center at 0 -3
\setquadratic \plot  2 3  1.5 1.5  0 0  -1.5 -1.5  -2 -3 /
        \plot  -2 3  -1.5 1.5  0 0  1.5 -1.5  2 -3 /
\put {$\bullet$} at 0 0
\put {$Q$} at 1 0
\putrectangle corners at 18 -3 and 28 3
\ellipticalarc axes ratio 3:1.1  360 degrees from 25 0 center at 23 0
\ellipticalarc axes ratio 3:1.1  360 degrees from 27 0 center at 23 0
\putrectangle corners at 8 -3 and 14 -7
\put {$\bullet$} at 11 -5
\put {$P$} at 11.8 -5
\arrow <.3truecm> [.2,.6] from 5 8 to 2 6
\arrow <.3truecm> [.2,.6]  from 17.5 6 to 20.5 4
\arrow <.3truecm> [.2,.6]  from 11 2.5 to 11 -1.5
\arrow <.3truecm> [.2,.6] from 4 -2 to 7 -5
\arrow <.3truecm> [.2,.6]  from 17 -2 to 15 -4
\put {$\pi_2$} at  3.5 7.9
\put {$h$} at  19.1 5.8
\put {$\pi$} at  11.7 .4
\put {$\pi_1$} at  5.7 -2.4
\put {$p$} at  16 -2
\put {$Y_2$} at  16.9 13
\put {$Y_1$} at -3.7 0
\put {$X$} at 14.7 -6
\put {$X^m$} at 29 0
\endpicture}
\vskip 1truecm
\centerline{\smc Figure 4}
\vskip 1truecm

\noindent
Denoting as usual
$$
K_{Y_2}=h^*(K_{X^m}+E_2)+(\nu_1-1)E_1+(\nu_2-1)E_2 \quad\text{and}\quad h^*(a_2E_2)= N_1E_1+N_2E_2
$$
we have clearly that $\nu_2=0$ and $N_2=-1$, and one computes that 
$\nu_1=\frac 15$ and $N_1=-\frac 15$. So
$$\align
z_P(X;s)&=\frac{\chi(C)}{(\nu_1+sN_1)(\nu_2+sN_2)} +
\frac{\chi(E_1\setminus C)}{\nu_1+sN_1}   
+\frac{\chi(E_2\setminus C)}{\nu_2+sN_2}  \\
&= \frac{-4}{(\frac15 - \frac15 s)(-s)} +
\frac{-4}{\frac15 - \frac15 s} + \frac 7{-s} = \frac {13}s \, ,
\endalign
$$
yielding $\lim_{s\to 1} z_P(X;s)= z_P(X;1)=13$.
Moreover on the level of the Grothendieck ring we have
$$
\align
\Cal Z_P(X;s)&=[C]\frac{(L-1)^2}{(L^{\nu_1+sN_1}-1)(L^{\nu_2+sN_2}-1)}
\\
&\qquad\qquad + L[C]\frac{L-1}{L^{\nu_1+sN_1}-1}   
+(L^2+L+1-[C])\frac{L-1}{L^{\nu_2+sN_2}-1}  \\
& = \frac {L^3-1+(L-1)[C]\sum_{i=1}^4 L^{\frac{i}{5}(1-s)}}{L^{-s}-1} 
\, ,
\endalign
$$
so \lq$\lim_{s\to 1} \Cal Z_P(X;s)$\rq\ $= \Cal Z_P(X;1)=-(L^3+L^2+L+4L[C])$.

\bigskip
\noindent
{\bf 3.6.1.} {\sl Note.} This last expression specializes to 
$$
-\big( (uv)^3+5(uv)^2-12u^2v-12uv^2 +5uv\big)
$$
on the Hodge polynomial level; as in [Ve3, 6.8 ] it is remarkable that al coefficients of the opposite polynomial have the \lq right\rq\ sign.

\bigskip
\noindent
{\bf 3.7.}  As promised we will verify that the negativity condition on the log discrepancies in (3.1) is equivalent to the absence of strictly lc singularities.  This is based on Lemma 3.7.1 below, which is an easy consequence of the following fact.
\bigskip
\proclaim
{3.7.0. Lemma} {\rm [KM, Lemma 3.39]} Let $p:W\rightarrow V$ be a proper birational morphism between normal varieties. Let $D$ be a $p$--nef $\Bbb Q$--Cartier divisor on $W$. Then 

(1) $-D$ is effective if and only if $p_*(-D)$ is.

(2) Assume that $-D$ is effective. Then for every $v \in V$ either $p^{-1}\{v\} \subset \supp D$ or $p^{-1}\{v\} \cap \supp D = \emptyset$.
\endproclaim
\bigskip
\proclaim
{3.7.1. Lemma}  Let $V$ be a normal variety and $B$ a $\Bbb Q$--divisor on $V$ such that $K_V + B$ is $\Bbb Q$--Cartier.  Let $p : V^m \rightarrow V$ be a log minimal model of $(V,B)$ as defined in (1.6.1(1)); we denote by $F = \sum_{i  \in T^m} F_i$ the reduced exceptional divisor of $p$, and by $B^m$ the strict transform of 
$B$ in $V^m$.  Also $a_i, i \in T_m$, is the log discrepancy of $F_i$ with respect to $(V,B)$; so
$$K_{V^m} + B^m + F = p^\ast (K_V + B) + \sum_{i \in T_m} a_i F_i.$$
Then
 (i) $a_i \leq 0$ for all $i \in T_m$, and

(ii) if $a_j < 0$, then $a_i < 0$ for all $F_i$ that satisfy $p(F_i) \subset p(F_j)$.
\endproclaim
\bigskip
\noindent
\demo{Proof}  (i) The divisor $D:=\sum_{i \in T_m} a_i F_i$ is $\Bbb Q$--Cartier and $p$--nef since $V^m$ is $\Bbb Q$--factorial and 
$K_{V^m}  + B^m + F$ is $p$--nef, respectively. Now $p_*(-D)=0$ and hence $-D$ is effective by Lemma 3.7.0(1).

(ii) Take a point $Q \in F_i \setminus \cup_{\ell \ne i} F_\ell$; then $p(Q) \in p(F_i) \subset p(F_j)$. Since $a_j<0$ we have by Lemma 3.7.0(2) that $p^{-1}\{p(Q)\} \subset \supp D$. In particular $Q \in \supp D$, implying that $a_i<0$.
\qed 
\enddemo

\bigskip
\noindent
{\bf 3.7.2.} {\sl Remark.}  The previous lemma is more generally valid for a $d$--minimal model $p : V^m_d \rightarrow V$ of $(V,B)$ as defined in (1.6.1(1)), where $d \in \Bbb Q$ and $0 \leq d \leq 1$.  But then the $a_i, i \in T$, are defined by
$$K_{V^m_d} + B^m + dF = p^\ast (K_V + B) + \sum_{i \in T_m} a_i F_i.$$
The proof is exactly the same since in this setting $K_{V^m_d} + B^m + dF$ is $p$--nef.
\bigskip
\proclaim
{\bf 3.7.3. Proposition}  Let $(X,B)$ be as in (3.1), and let $p : X^m \rightarrow X$ be a log minimal model of $(X,B)$.  Then $(X,B)$ has no strictly lc singularities if and only if all log discrepancies with respect to $(X,B)$ of exceptional divisors of $p$ are (strictly) negative.
\endproclaim
\bigskip
\noindent
\demo{Proof}  We use the notations $F = \sum_{i \in T^m} F_i$ and $B^m$ from (3.1).  The log discrepancies $a_i$, $i \in T$, are given by
$$K_{X^m} + B^m + F = p^\ast (K_X + B) + \sum_{i \in T_m} a_iF_i.$$ 
We will show that $(X,B)$ has a strictly lc singularity if and only if some $a_i, i \in T_m$, is zero.  (Since all $a_i \leq 0$ by Lemma 3.7.1(i), this is clearly equivalent to the assertion.)

Suppose first that $Q \in X$ is a strictly lc singularity of $(X,B)$.  Then $Q$ has a neighbourhood $U$ such that $(U,B|_U)$ is lc.  Hence all $F_\ell$ occurring on $p^{-1}U$ must have $a_\ell = 0$.  Furthermore there is a least one such $F_\ell$, because otherwise $p|_{p^{-1}U}$ would have no exceptional divisors, and this would mean that $(U,B|_U)$ is klt.

On the other hand, suppose that $a_i = 0$ for some $i \in T_m$.  Then, by Lemma 3.7.1(ii), $a_j = 0$ for all $F_j$ that satisfy 
$p(F_i) \subset p(F_j)$.  This implies that $\cup_{\Sb \ell \in T^m \\ a_\ell < 0 \endSb} p(F_\ell)$ intersects $p(F_i)$ in a proper closed subset.  Take now any $Q$ in $p(F_i) \setminus \cup_{\Sb \ell \in T^m \\ a_\ell < 0 \endSb} p(F_\ell)$.  Then $Q$ has a neighbourhood $U$ (in $X$) such that only $F_\ell$ with $a_\ell = 0$ occur in $p^{-1}U$.  This easily yields that $(U,B|_U)$ is lc.  And since $a_i = 0$ we have that $Q$ is a strictly lc singularity of $(X,B)$.  \qed
\enddemo
\bigskip
\bigskip
\head
4. Stringy zeta functions for arbitrary pairs
\endhead
\bigskip
\noindent
{\bf 4.1.} Let again $(X,B)$ be as in (3.1).  In the special case that $(X,B)$ has some strictly lc singularity we cannot associate a stringy zeta function to $(X,B)$ as in (3.2); see (2.6).  A solution is to use, instead of a (usual) log minimal model $X^m$ of $X$, a $d$--minimal model $X^m_d$ as defined in (1.6.1) for some $d < 1$.  Alternatively, we can use here as well the $d$--canonical model $X^c_d$, which has the advantage of being unique.  We will work with $X^c_d$; see (4.5) with $X^m_d$.
\bigskip
\noindent
{\bf 4.2. Definition.}  We assume the (relative log) MMP.  Fix $d \in \Bbb Q$ with $0 \leq d < 1$.  To any pair $(X,B)$ as in (3.1) we associate the stringy zeta function
$$\Cal Z_d(s) = \Cal Z_d((X,B);s) := \Cal Z((X^c_d, B^c + dF^c), K_{X^c_d} + B^c + dF^c - q^\ast (K_X + B); s).$$
Here $q : X^c \rightarrow X$ is the (relative) $d$--canonical model of $(X,B)$ as defined in (1.6.1 (2)) with reduced exceptional divisor $F^c$, $B^c$ is the strict transform of $B$ by $q$, and the right hand side is the zeta function associated in (2.2) to the klt pair $(X^c_d, B^c + dF^c$) and the \lq log discrepancy divisor\rq\ $K_{X^c_d} + B^c + dF^c - q^\ast (K_X + B)$ on $X^c_d$. (When $B=0$ we just write $\Cal Z_d(X;s)$.)
\bigskip
\noindent
{\bf 4.2.1. Formula.}  Take a log resolution $g : Y \rightarrow X^c$ of the pair $(X^c_d, \supp B \cup \supp F^c)$.  Let $E_i, i \in T$, be the irreducible components of the exceptional divisor of $g$ and of the strict transform of $B^c + F^c$, and denote as usual $E^\circ_I = (\cap_{i \in I} E_i) \setminus (\cup_{\ell \notin I} E_\ell)$ for $I \subset T$.  Then
$$\Cal Z_d(s) = \sum_{I \subset T} [E^\circ_I] \prod_{i \in I} \frac{L-1}{L^{\nu_i + sN_i} - 1} ,$$
where
$$\align
&K_Y = g^\ast(K_{X^c_d} + B^c + dF^c) + \sum_{i \in T} (\nu_i - 1) E_i \qquad\qquad \text{and} \\
&g^\ast (K_{X^c_d} + B^c + dF^c - q^\ast (K_X + B)) = \sum_{i \in T} N_i E_i.
\endalign$$
Remark that here the $\nu_i (\in \Bbb Q_{> 0})$ and the $N_i( \in \Bbb Q)$ depend on $d$.    
\bigskip
\noindent
{\bf 4.2.2.} (i) Of course analogously we can associate to $(X,B)$ stringy zeta functions on other levels, e.g. $z_d (s)$ with Euler characteristics. 

\noindent
(ii) For a constructible subset $W$ of $X$ we can introduce more generally $\Cal Z_{d,W}(s) = \Cal Z_{d,W} ((X,B);s)$ and for instance $z_{d,W}(s)$, analogously as in 3.2.3(ii).
\bigskip
\noindent
{\bf 4.3.}  When $B=0$ and $d=0$ the variety $X^c_0$ is just the relative canonical model of $X$ and $\Cal Z_0(s) = \Cal Z_0(X;s) = \Cal Z((X^c_0,0), K_{X^c_0} - q^\ast K_X;s)$.

In the context of generalizing the elliptic genus to singular varieties, Totaro [To] also used the relative canonical model in a similar way.
\bigskip
\noindent
{\bf 4.4.}  Again we make a link with question (I) in the introduction.  Consider a pair $(X,B)$ as in (3.1) that has a log resolution $\pi : Y \rightarrow X$ for which {\sl all} log discrepancies with respect to $(X,B)$ of divisors on $Y$ are nonzero; i.e.
$$K_Y =  \pi^\ast (K_X + B) + \sum_{i \in T} (a_i - 1)E_i \quad \text { with all } a_i \ne 0,$$
where the $E_i, i \in T$, are as usual.  Take $d \in \Bbb Q$ with $0 \leq d < 1$ and assume that $\pi$ factorizes through the $d$--canonical model $q : X^c_d \rightarrow X$.  Denote $\nu_i, N_i$ as in (4.2.1).  By the same computation as in (3.4.1) we have that $\nu_i + N_i = a_i$ (which does not depend on $d$), and
$$\Cal Z_d(1) = \sum_{I \subset T}  [E^\circ_I] \prod_{i \in I} \frac{L-1}{L^{a_i}-1} ,$$
which is again the `invariant' of (3.4.2).
\bigskip
\noindent
{\bf 4.5.}  In definition 4.2 we could have taken any (relative) $d$--minimal model $p : X^m_d \rightarrow X$ of $(X,B)$ instead of $X^c_d$.  
\bigskip
\proclaim
{Proposition}  For $(X,B)$ as in (3.1), let $p : X^m_d \rightarrow X$ be a $d$--minimal model of $(X,B)$.  Let $F^m$ and $B^m$ denote the reduced exceptional divisor of $p$ and the strict transform of $B$ by $p$, respectively.  Then
$$\Cal Z_d(s) = \Cal Z((X^m_d,B^m + dF^m), K_{X^m_d} + B^m +dF^m - p^\ast(K_X + B); s),$$
where the right hand side is the zeta function associated in (2.2) to the klt pair $(X^m_d,B^m + dF^m)$ and the divisor $K_{X^m_d} + B^m + dF^m - p^\ast (K_X + B)$ on $X^m_d$. 
\endproclaim
\bigskip
\demo{Proof}                
We still use the notation of (4.2).  Consider the diagram $p : X^m_d \overset f \to \rightarrow X^c_d \overset q \to \rightarrow X$ (see (1.6.3)).  We claim that
$$K_{X^m_d} + B^m + dF^m = f^\ast (K_{X^c_d} + B^c + dF^c).$$
Indeed, let $\sum_{i \in S} F_i$ be the reduced exceptional divisor of $f$ and put
$$K_{X^m_d} + B^m + dF^m = f^\ast (K_{X^c_d} + B^c + dF^c) + \sum_{i \in S} a_i F_i.$$
First remark that $f: X^m_d \rightarrow X^c_d$ is a $d$--minimal model of $(X^c_d, B^c + dF^c)$.
Then all $a_i \leq 0$ by Lemma 3.6.1(i) and Remark 3.6.2.  On the other hand, also $a_i \geq 0$ since $\logdisc (X^c_d, B^c + dF^c) \geq 1 - d$ (by definition of $X^c_d$).

\smallskip
Take now a log resolution $h : Y \rightarrow X^m_d$ of $(X^m_d, B^m + dF^m)$ such that $g := f \circ h$ is also a log resolution of $(X^c_d, B^c + dF^c)$.  Then, with the notation of (4.2.1), we have that
$$K_Y = g^\ast (K_{X^c_d} + B^c + dF^c) + \sum_{i \in T} (\nu_i - 1) E_i = h^\ast (K_{X^m_d} + B^m + dF^m) + \sum_{i \in T} (\nu_i - 1)E_i$$
and
$$h^\ast (K_{X^m_d} + B^m + dF^m - p^\ast (K_X + B)) = g^\ast (K_{X^c_d} + B^c + dF^c - q^\ast (K_X + B)) = \sum_{i \in T} N_i E_i.$$
This proves the assertion. \qed
\enddemo
\bigskip
\noindent
Note that in particular the right hand side in the proposition does not depend on the chosen $d$--minimal model $X^m_d$; alternatively this can be verified as in (3.2.2).
\bigskip
\noindent
{\bf 4.6.} 
When the pair $(X,B)$ has no strictly lc singularities, we can associate to it both the zeta function $\Cal Z(s)$ of (3.2) and, for $0 \leq d < 1$, the zeta functions $\Cal Z_d(s)$ of (4.2).  Some natural questions arise in this context.  We first give an example.

\vskip 1true cm
\centerline{\beginpicture
\setcoordinatesystem units <.5truecm,.5truecm>
\putrectangle corners at 9 11 and 18 6
\putrule from 10 8 to 16 8
\putrule from 11 7 to 11 10
\putrule from 13 7 to 13 10
\putrule from 15 7 to 15 10
\multiput {$E_2$} at 16.7 8  25.9 2 /
\multiput {$E_1$} at 10.3 9.5   1.2 3.5  /
\multiput {$E_0$} at 14.3 9.5   6.8 3.5   23.5 3.5   11.4 -5.2 /
\multiput {$B$} at 12.4 9.5  4.6 4  21.7 3.5  14.2 -2.3  12.3 -10.2 /
\putrectangle corners at 0 0 and 8 5
\putrule from 4 1 to 4 4
\setlinear  \plot 2 1.5  6 3.5 /
            \plot 6 1.5  2 3.5 /
\putrectangle corners at 19 0 and 27 5
\putrule from 22.3 1 to 22.3 4
\putrule from 24.2 1 to 24.2 4
\putrule from 19.9 2 to 25.2 2
\multiput {$\bullet$} at 20.9 2  13.5 -11.5 /
\put {$Q$} at 20.7 1.3
\putrectangle corners at 10 -1 and 17 -6
\putrule from 11 -4.5 to  16 -4.5
\setquadratic \plot 12 -2   13.5 -4.5   15 -2 /  
\putrectangle corners at 10.5 -9 and 16.5 -14
\ellipticalarc axes ratio 4:3  90 degrees from 11.5 -10 center at 13.5 -10
\ellipticalarc axes ratio 4:3  90 degrees from 15.5 -13 center at 13.5 -13
\put {$P$} at 13.5 -12.2
\arrow <.3truecm> [.2,.6] from 8 8 to 5 6
\arrow <.3truecm> [.2,.6]  from 19 8 to 22 6
\arrow <.3truecm> [.2,.6]  from 13.5 5 to 13.5 0
\arrow <.3truecm> [.2,.6] from 6 -1 to 9 -3
\arrow <.3truecm> [.2,.6]  from 21 -1 to 18 -3
\arrow <.3truecm> [.2,.6]  from 13.5 -6.6 to 13.5 -8.4
\put {$f_2$} at  6.6 8
\put {$c_2$} at  20.4 7.8
\put {$f$} at  14.1 2.5
\put {$f_1$} at  6.7 -2.4
\put {$c_1$} at  20.2 -2.4
\put {$f_0$} at  14.1 -7.2
\put {$Y_2$} at  18.7 10
\put {$Y_1$} at -.7 2.5
\put {$Y_0$} at 17.7 -5
\put {$\tilde Y$} at 27.7 2.5
\put {$X$} at 17.2 -11.5
\endpicture}
\vskip 1truecm
\centerline{\smc Figure 5}
\vskip 1truecm

\noindent
{\bf 4.6.1.} {\sl Example.}  Let $P \in X$ be a simple elliptic surface singularity (germ), i.e. the exceptional divisor of its minimal resolution $f_0 : Y_0 \rightarrow X$ is just one (smooth) elliptic curve $E_0$ with self--intersection number $-\kappa_0$ on $Y_0$.  Let $B \ni P$ be an irreducible divisor on $X$ whose strict transform by $f_0$ intersects $E_0$ in just one point with intersection multiplicity $2$; see Figure 5.
Further we will use the same notation for curves and their strict transforms.

Let $f : Y_2 \overset f_2 \to \longrightarrow Y_1 \overset f_1 \to \longrightarrow Y_0$ be the minimal log resolution of $(Y_0, E_0 \cup B)$.  Here $f_i$ is a blowing--up with exceptional curve $E_i$.  The log discrepancies $a_i$ of $E_i$ with respect to $(X, \frac 12 B)$ are easily computed as $a_0 = - \frac{1}{\kappa_0}$, $a_1 = - \frac{1}{\kappa_0} + \frac 12$ and $a_2 = - \frac{2}{\kappa_0}$.  So $a_0 < 0$ and $a_2 < 0$; on the other hand $a_1$ can be negative or positive and in one case zero (when $\kappa_0 = 2$).  In particular $P \in (X, \frac 12 B)$ is not log canonical.

Since $E_1$ has self--intersection number $-2$ on $Y_2$ we can consider the contraction $c_2 : Y_2 \rightarrow \tilde Y$, mapping $E_1$ to the $A_1$--singularity $Q \in \tilde Y$.  Finally $c_1 : \tilde Y \rightarrow Y_0$ is then the contraction of $E_2$.
\bigskip
\noindent
{\bf Claim 1 :} {\sl $f_0 \circ c_1 : \tilde Y \rightarrow X$ is a log minimal model of $(X, \frac 12 B)$.}
\medskip
\noindent
Indeed, $(\tilde Y, \frac 12 B + E_0 + E_2)$ is clearly dlt.  We check that
$(K_{\tilde Y} + \frac 12 B + E_0 + E_2) \cdot E_i \geq 0$ for $i = 0,2$ :
$$(K_{\tilde Y} + \frac 12 B + E_0 + E_2) \cdot E_0 =  \deg K_{E_0} + \frac 12 B \cdot E_0 + E_2 \cdot E_0 = 0 + 0 + 1 = 1,$$
and
$$(K_{\tilde Y} + \frac 12 B + E_0 + E_2) \cdot E_2 = \deg(K_{E_2}
 + \Diff) + \frac 12 B \cdot E_2 + E_0 \cdot E_2 = (-2 + \frac 12) + \frac 12 + 1 = 0.$$
Here we used the notion of {\it Different}, see e.g. [Kol, 16.5-6].  Alternatively one can use that $c^\ast_2 E_2 = E_2 + \frac 12 E_1$ and $c^\ast_2 K_{\tilde Y} = K_{Y_2}$.

\smallskip
We now compute the stringy zeta function of $(X, \frac 12 B)$, or rather of its germ in $P$; for simplicity we work on the level of Euler characteristics $z_P(s)$.  We take $X^m = \tilde Y$ in Definition 3.2.  Then
$$z_P(s) = z_P ((\tilde Y, \frac 12 B + E_0 + E_2), K_{\tilde Y} + \frac 12 B + E_0 + E_2 - (f_0 \circ c_1)^\ast (K_X + \frac 12 B);s).$$
We take $c_2$ as the resolution $h$ in Formula 3.2.1.  For the numbers $\nu_i$ and $N_i$ in this formula we have clearly that
$$\alignat 3
\nu_0 &= 0, &\qquad \nu_2 &= 0, &\qquad \nu_B &= \frac12, \\N_0 &= a_0 = -\frac{1}{\kappa_0}, &\qquad N_2 &= a_2 = - \frac{2}{\kappa_0},&\qquad  N_B &= 0,
\endalignat
$$
and then it is easy to verify that $\nu_1 = \frac 12$ and $N_1 = -\frac {1}{\kappa_0}$.  Hence
$$\align
z_P(s) & = \frac{1}{\nu_2 + sN_2} \big(-1 + \frac{1}{\nu_1 + sN_1} + \frac{1}{\nu_B + sN_B} + \frac{1}{\nu_0 + sN_0}\big) + \frac{1}{\nu_1 + sN_1} + \frac{-1}{\nu_0 + sN_0} \\
& = \frac{\kappa_0^2}{2s^2} - \frac{\kappa_0}{2s}.
\endalign
$$
(Alternatively one can use the shorter formula in [Ve1, \S4], which also applies to our stringy zeta functions.)
\bigskip
\noindent
{\bf Claim 2 :} {\sl $f_0 : Y_0 \rightarrow X$ is both $d$--minimal model and $d$--canonical model of $(X, \frac 12 B)$ for $0 \leq d < 1$.}
\medskip
\noindent
One easily verifies that
$$K_{Y_2} + \frac 12 B + dE_0 + E_1 + E_2 = f^\ast (K_{Y_0} + \frac 12 B + dE_0) + (\frac 32 - d)E_1 + 2(1 - d)E_2. \tag $*$ $$
This implies that $\logdisc (Y_0, \frac 12 B + dE_0) > 1 - d$.  And $K_{Y_0} + \frac 12 B + dE_0$ is also $f_0$--ample since
$$(K_{Y_0} + \frac 12 B + dE_0) \cdot E_0 = \frac 12 B \cdot E_0 - (1-d) E^2_0 = 1 + \kappa_0 (1-d) > 1 \, (> 0).$$
We take $f$ as the log resolution $g$ in Formula 4.2.1.  Here for the numbers $\nu_i$ and $N_i$ in this formula we obtain, using ($\ast$) and the fact that $a_i = \nu_i + N_i$ :
$$\alignat4
 \nu_B &= \frac 12, &\qquad \nu_0 &= 1-d, &\qquad  \nu_1 &= \frac 32 - d, &\qquad  \nu_2 &= 2(1-d), \\
  N_B &=  0,  &\qquad   N_0 &= d-1-\frac{1}{\kappa_0}, 
  &\qquad  N_1 &= d-1- \frac{1}{\kappa_0}, 
  &\qquad  N_2 &= 2(d - 1 - \frac{1}{\kappa_0}).
\endalignat
$$
Hence
$$\align
z_{d,P}(s) & = \frac{1}{\nu_2 + sN_2} (-1 + \frac{1}{\nu_1 + sN_1} + \frac{1}{\nu_B + sN_B} + \frac{1}{\nu_0 + sN_0}) + \frac{1}{\nu_1 + sN_1} + \frac{-1}{\nu_0 + sN_0} \\
& = \frac{1}{2(1 - d + (d - 1 - \frac{1}{\kappa_0})s)^2} + \frac{1}{2(1-d+(d-1 - \frac{1}{\kappa_0})s)} \, .
\endalign
$$
So in fact $z_P(s) = \lim_{d \rightarrow 1} z_{d,P} (s)$.  (This is even true for the $\nu_i$ and $N_i$.)
\bigskip
\noindent
{\bf 4.6.2.} (i)  In the previous example the log minimal model $X^m = \tilde Y$ is {\sl not} a $d$--minimal model for $d < 1$, not even for $d$ close to $1$.  At least for surfaces it is not obvious to give such examples; for instance when $B=0$ we verified that then a log minimal model $X^m$ is also a $d$--minimal model for $d$ close to $1$.

(ii)  In general, when $(X,B)$ has no strictly lc singularities, if a log minimal model $X^m$ is also a $d$--minimal model for $d$ close to $1$, one can check that $z(s) = \lim_{d \rightarrow 1} z_d(s)$.  It is remarkable that this is still true in Example 4.6.1.
\bigskip
\noindent
{\bf 4.6.3. Questions.}  (i) When is $X^m = X^m_d$ for $d$ close to $1$ ?

(ii) Is $z(s) = \lim_{d \rightarrow 1} z_d(s)$ ?
\bigskip
\noindent
{\bf 4.7.}  By Theorem 3.6.3, the only singularities which were not covered by Definition 3.2 are the strictly lc singularities.  We will determine the stringy zeta functions of Definition 4.2 for the strictly lc singularities $P \in X$ on a normal surface $X$.  Recall first their classification  [Al], given in Figure 6 by the dual graph of the minimal log resolution $\pi : Y \rightarrow X$ of $P \in X$.  The exceptional components $E_i, i \in T$, are represented by dots and an intersection between them by a line connecting the corresponding dots.  All components are rational, except in (1); and
in (4) the $n_i$ are the possible absolute values of the determinants of the intersection matrices of the three disjoint chains, separated by the central component $E$.
\vskip 1true cm
\centerline{\beginpicture
\setcoordinatesystem units <.5truecm,.5truecm>
\put {$\bullet$} at  3  18
\put {$(1)$} at  -5 18
\put{elliptic curve} at 16 18
\multiput {$\bullet$} at 3 10  5.1 10.9  6 13  5.1 15.1  3 16  
                         0.9 15.1  0 13  0.9 10.9   /
\plot 3 10  5.1 10.9  6 13  5.1 15.1  3 16  
                         0.9 15.1  0 13  0.9 10.9  3 10  /
\put {$(2)$} at  -5 13
\put{a closed chain of length $r\geq 2$} at 19 13
\multiput {$\bullet$} at  -2 7  -2 5  0 6  2 6  6 6  8 6  10 7  10 5 /
\putrule from 0 6 to 3 6 
\putrule from 5 6 to 8 6
\put {$\dots$} at 4 6
\plot -2 7  0 6 /           
\plot -2 5  0 6 /
\plot 10 7  8 6 /           
\plot 10 5  8 6 /
\put {$E_2$} at -2 4.3
\put {$E_1$} at -2 7.6
\put {$E_4$} at 10 4.3
\put {$E_3$} at 10 7.6
\put {$E_5$} at .2 5.3
\put {$E_6$} at 2.1 5.3
\put {$E_{5+k}$} at 7.8 5.3
\put {$(3)$} at  -5 6	  
\put {where $k \geq 0$, and $E_1,E_2,E_3,E_4$}	 
at 19.5 6.5
\put {have self--intersection $-2$}	 
at 18.2 5.5
\put {$E$} at 4 1.7
\put {$(4)$} at  -5 1
\multiput {$\bullet$} at 4 1  2.4 1  -.8 1  -2.4 1 
                              5.6 2  8.8 2  10.4 2
							  5.6 0  8.8 0  10.4 0 /
\putrule from -2.4 1 to 0 1
\putrule from 1.6 1 to 4 1
\putrule from 5.6 2 to 6.4 2
\putrule from 8 2 to 10.4 2
\putrule from 5.6 0 to 6.4 0
\putrule from 8 0 to 10.4 0
\multiput {$\dots$} at .8 1  7.2 2  7.2 0 /
\plot 5.6 2  4 1  5.6 0 /
\put{with $(n_1,n_2,n_3) = \left\{ \aligned & (2,3,6) \\ & (2,4,4) \\ & (3,3,3)  \endaligned \right.$}
at 19 1
\endpicture}
\vskip 1truecm
\centerline{\smc Figure 6}
\vskip 1truecm

\noindent
{\bf 4.8.}  Let $d \in \Bbb Q$ with $0 \leq d < 1$.  We compute the stringy zeta functions $\Cal Z_{d,P}(s) = \Cal Z_{d,P} (X;s)$ for the strictly lc surface singularities $P \in X$ as described in (4.7).
\bigskip
\noindent
{\bf Case (1).}  We have $K_Y = \pi^\ast K_X - E$.  Hence $(K_Y + dE) \cdot E = (1 - d)(-E^2) > 0$, meaning that $K_Y + dE$ is $\pi$--ample.  Also it is clear that $\logdisc (Y,dE) > 1-d$.  So in fact $Y = X^m_d = X^c_d$.  Then
$$\Cal Z_{d,P}(s) = [E] \frac{L-1}{L^{(1-d)+(d-1)s}-1} = [E] \frac{L-1}{L^{(1-d)(1-s)}-1} \, .$$
Since $\chi(E) = 0$ we have $z_{d,P}(s) = 0$. 
\bigskip
\noindent
{\bf Case (2).}  We have $K_Y = \pi^\ast K_X - \sum^r_{i=1} E_i$.  Hence
$$\align
(K_Y + d \sum^r_{i=1} E_i) \cdot E_j & =  (K_Y + \sum^r_{i=1} E_i) \cdot E_j - (1-d) \sum^r_{i=1} E_i \cdot E_j \\
& = 0 - (1-d) (E^2_j + 2) \geq 0
\endalign
$$
for each $j \in \{ 1, \cdots, r \}$. (Indeed, $E^2_j \leq -2$ because $\pi$ is the {\it minimal} log resolution of $P \in X$.)  So 
$K_Y + d \sum^r_{i=1} E_i$ is $\pi$--nef, and since also $\logdisc (Y,d \sum^r_{i=1} E_i) > 1-d$, we have $Y = X^m_d$.  Then by (4.5) we have
$$\align
\Cal Z_{d,P}(s) & = r(L-1) \frac{L-1}{L^{(1-d)+(d-1)s}-1} + r \frac{(L-1)^2}{(L^{(1-d)+(d-1)s}-1)^2} \\
& = r(L-1)^2 \Big(\frac{1}{(L^{(1-d)(1-s)}-1)^2} + \frac{1}{L^{(1-d)(1-s)}-1}\Big).
\endalign
$$ 

\noindent
{\bf Case (3).}  We have $K_Y = \pi^\ast K_X - \frac 12 \sum^4_{i=1} E_i - \sum^{5+k}_{i=5} E_i$.  For $0 < d < 1$ we claim that $X^m_d$ is obtained from $Y$ by contracting $E_1,E_2,E_3,E_4$.  Denote for the moment this contraction by $h : Y \rightarrow X^\prime$ and $F_i := h(E_i)$ for $i = 5, \cdots , 5+k$.  So $\pi$ factorizes as $Y \overset h \to \longrightarrow X^\prime \overset p \to \longrightarrow X$.  We first treat the general case $k > 0$.
\medskip
($i$) $\logdisc (X^\prime, d \sum^{5+k}_{i=5} F_i) > 1-d$.  
\medskip
\noindent
It is easy to see that $h^\ast F_5 = E_5 + \frac 12 E_1 + \frac 12 E_2$ and $h^\ast F_{5+k} = E_{5+k} + \frac 12 E_3 + \frac 12 E_4$, and that $K_Y = h^\ast K_{X'}$.  Hence
$$K_Y + d \sum^{k+5}_{i=5} E_i - h^\ast (K_{X^\prime} + d \sum^{5+k}_{i=5} F_i)  = d \sum^{5+k}_{i=5} (E_i - h^\ast F_i)   = \sum^4_{i=1} - \frac d2 E_i .
$$
So the log discrepancies of 
$E_i, 1  \leq i \leq 4$, are $1 - \frac d2 > 1-d$.
\medskip
($ii$) $K_{X^\prime} + d \sum^{5+k}_{i=5} F_i$ is $p$--nef.
\medskip
\noindent
We have that
$$\align (K_{X^\prime} + d \sum^{5+k}_{i=5} F_i) \cdot F_5 & = h^\ast (K_{X^\prime} + d \sum^{5+k}_{i=5} F_i) \cdot E_5 \\
& = (K_Y + d \sum^{5+k}_{i=5} E_i + \frac d2 \sum^4_{i=1} E_i) \cdot E_5 \\
& = -2 - E^2_5 + d(E^2_5 + 1) + \frac d2 \cdot 2 \\
& = (1-d) (-E^2_5 - 2) \geq 0.
\endalign
$$
Analogously $(K_{X^\prime} + d \sum^{5+k}_{i=5} F_i)  \cdot F_\ell \geq 0$ for $5 < \ell \leq 5+k$.  

\smallskip
So indeed $X^\prime = X^m_d$.  Then, by (4.5) and with the usual $\nu_i$ and $N_i$, we have 
$$\align z_{d,P}(s) = & \sum^{4+k}_{i=5} \frac{1}{(\nu_i + sN_i) (\nu_{i+1} + sN_{i+1})} + \frac{-1}{\nu_5 + sN_5} + \frac{-1}{\nu_{5+k} + sN_{5+k}} \\
& + \sum^2_{i=1} \frac{1}{\nu_i + sN_i} (1 + \frac{1}{\nu_5 + sN_5}) + \sum^4_{i=3} \frac{1}{\nu_i + sN_i} (1 + \frac{1}{\nu_{5+k} + sN_{5+k}}).
\endalign
$$      
For $1 \leq i \leq 4$ we already saw that $\nu_i = 1 - \frac d2$, and it is easy to check that $N_i = \frac{d-1}2$.  Hence
$$\align
z_{d,P}(s)  = & k \frac{1}{((1-d)+(d-1)s)^2} + 2 \frac{-1}{1-d+(d-1)s} \\ & + 4 \frac{1}{(1 - \frac d2 + (\frac{d-1}2)s)} \Big( 1 + \frac{1}{(1-d)+(d-1)s)} \Big) \\
 = & \frac{k}{(1-d)^2(1-s)^2} + \frac{6}{(1-d)(1-s)}.
\endalign
$$  
A slightly longer computation yields
$$\Cal Z_{d,P}(s) = \frac{k(L-1)^2}{(L^{(1-d)(1-s)}-1)^2} + \frac{(L-1)\bigl((k-1)(L-1) + 2L + 4L^{1 - \frac d2 + \frac{d-1}2 s}\bigr)}{L^{(1-d)(1-s)}-1}.$$
For the case $k=0$ we obtain analogously
$$z_{d,P} (s) = \frac{6}{(1-d)(1-s)}$$
and
$$\Cal Z_{d,P}(s) = \frac{(L-1)(L+1+4L^{1 - \frac d2 + \frac{d-1}2 s})}{L^{(1-d)(1-s)} -1}$$
which is compatible with the formulas for $k > 0$.

When $d=0$ the ordinary (relative) minimal model $X^m_0$ of $P \in X$ is $Y$, but one computes that $z_{0,P}(s)$ and $\Cal Z_{0,P}(s)$ are just given by putting $d=0$ in the previous formulas.
\bigskip
\noindent
{\bf Case (4).}  Let $h : Y \rightarrow X^\prime$ be the contraction of all components $E_i, i \in T$, except $E$, and denote $F := h(E)$.  So $\pi$ factorizes as $Y \overset h \to \longrightarrow X^\prime \overset p \to \longrightarrow X$.

We have that $K_{X^\prime} + dF$ is always $p$--ample.  Indeed :
$$\align (K_{X^\prime} + dF) \cdot F & = (K_{X^\prime} + F) \cdot F - (1 - d)F^2 \\ & = -2 + \deg(\Diff) - (1-d)F^2 \\ & = -2 + \sum^3_{i=1} \frac{n_i - 1}{n_i} - (1-d)F^2 = 0 + (1-d)(-F^2) > 0,
\endalign
$$
using again the Different (see [Kol, 16.5-6]).

On the other hand, the concrete expressions for the divisors $K_Y - h^\ast K_{X^\prime}$ and $h^\ast F$ depend on the self--intersection numbers of the $E_i, i \in T$, but one can verify that in each case $\logdisc (X^\prime, dF) \geq 1-d$, for $d$ close enough to 1.  So for such $d$ we have that $X^\prime = X^c_d$, and then by (4.2.1) and the formula in [Ve1, \S4] (which is also valid in this context) we obtain 
$$z_{d,P}(s) = \frac{1}{1-d+(d-1)s} (-1 + n_1 + n_2 + n_3) = \frac{n_1  + n_2 + n_3 - 1}{(1-d)(1-s)}.$$
The expression for $\Cal Z_{d,P}(s)$ depends on the concrete case, but can be given simultaneously for all cases (with still $d$ close enough to 1) in terms of the determinants $D_i$ of [Ve1, \S5]:
$$\Cal Z_{d,P}(s) = \frac{L-1}{L^{(1-d)(1-s)}-1} (L - 2 + D_1 + D_2 + D_3).$$

\vskip 1true cm
\centerline{\beginpicture
\setcoordinatesystem units <.5truecm,.5truecm>
\putrule from -2 0 to 2 0
\putrule from 0 0 to 0 -2 
\multiput {$\bullet$} at  0 0  2 0  -2 0  0 -2 /
\put {$E$} at  0 .7
\put {$E_1$} at -2.7 0
\put {$E_3$} at 2.7 0
\put {$E_2$} at -.5 -2.4
\endpicture}
\vskip 1truecm
\centerline{\smc Figure 7}
\vskip 1truecm

We illustrate this with the concrete example of Figure 7,
where $E^2_1 = E^2_2 = E^2_3 = -3$.   
Here $K_Y = h^\ast K_{X^\prime} - \frac 13 \sum^3_{i=1} E_i$ and $h^\ast F = E + \frac 13 \sum^3_{i=1} E_i$, yielding
$$K_Y + dE - h^\ast(K_{X^\prime} + dF) = \sum^3_{i=1} \Big(\frac{2-d}3 - 1\Big)E_i.$$
So $\logdisc (X^\prime, dF) \geq 1-d$ if and only if $d \geq \frac 12$, and then
$$z_{d,P}(s) = \frac{8}{(1-d)(1-s)} \quad  \text { and } \quad \Cal Z_{d,P}(s) = \frac{L-1}{L^{(1-d)(1-s)} - 1} (L - 2 + 3D),$$
where 
$$D = 1 + L^{\frac{2-d}{3} + \frac{d-1}3 s}+ (L^{\frac{2-d}3 + \frac{d-1}3 s})^2.$$
(This can easily be verified without using [Ve1].)

When $0 \leq d \leq \frac 12$ then clearly $\logdisc (Y,dE + \sum^3_{i=1} dE_i)> 1-d$, and moreover $K_Y + dE + \sum^3_{i=1} dE_i$ is $\pi$--nef since
$$(K_Y + dE + d \sum^3_{i=1} E_i) \cdot E = -2 - (1-d)E^2 + 3d \geq d \geq 0$$
and
$$(K_Y + dE + d \sum^3_{i=1} E_i) \cdot E_j = 1 - 2d \geq 0$$
for $j = 1,2,3$. (We used the fact that $E^2 \leq -2$.)  So now $Y = X^m_d$ and then by (4.5) we obtain
$$\align
z_{d,P}(s) & = \frac{-1}{1-d+(d-1)s} + 3 \frac{1}{1-d+(d-\frac 23)s}\\ &\qquad\qquad\qquad + 3 \frac{1}{(1-d+(d-1)s)(1-d+(d-\frac 23)s)} \\
& = \frac{5-2d+(2d - \frac 73)s}{(1-d)(1-s)(1-d+(d- \frac 23)s)}
\endalign
$$
and
$$\align \Cal Z_{d,P}(s)  = & (L-2) \frac{L-1}{L^{(1-d)+(d-1)s}-1} + 3L \frac{L-1}{L^{1-d+(d-\frac 23)s} - 1} \\ & + 3 \frac{(L-1)^2}{(L^{(1-d)+(d-1)s}-1)(L^{1-d+(d-\frac 23)s}-1)} \\
 = & \frac{(L-1)(-1-L+3L^{2-d+(d-1)s} + (L-2) L^{1-d+(d- \frac 23)s} )}{(L^{(1-d)(1-s)}-1)(L^{1-d+(d - \frac 23)s} - 1)} .
\endalign
$$        
For $d = \frac 12$ these expressions are indeed the same as the previous ones.
\bigskip
\bigskip
\head
5. Stringy invariants without MMP
\endhead
\bigskip 
\noindent
{\bf 5.1.} Here we present a partial result concerning question (I).  Let $X$ be a quasiprojective $\Bbb Q$--Gorenstein variety.  With the notation of (0.1), we associate to a log resolution $\pi : Y \rightarrow X$ of $X$ with {\it all\/} log discrepancies 
$a_i \ne 0$ the `stringy expression'
$$\sum_{I \subset T} [E^\circ_I] \prod_{i \in I} \frac{L-1}{L^{a_i}-1}.$$
Restricting to log resolutions $\pi$ {\it that factorize through the blowing--up of $X$ in $X_{\sing}$}, we will show that this expression is indeed an invariant of $X$.
\bigskip
\noindent
\proclaim{5.2. Lemma}  Let $X$ be a quasiprojective normal variety and $p : \tilde X \rightarrow X$ the blowing--up of $X$ in $X_{\sing}$.  Then there exists a linear system $\Cal L$ on $X$ with base locus $X_{\sing}$, such that the induced linear system $p^\ast \Cal L - E$ on $\tilde X$ is base--point free, where $E$ is the exceptional divisor of $p$.
\endproclaim
\bigskip
\demo{Proof}  Suppose first that $X$ is projective and take a very ample sheaf $\Cal M$ on $X$.  Denoting by $\Cal I$ the ideal sheaf of $X_{\sing}$ in $X$, we have that $\Cal I \otimes \Cal M^t$ is generated by global sections for $t$ large enough, see e.g. [Ha, Theorem II 5.17].  Then by [BS, Theorem 2.1] the sheaf $p^\ast (\Cal I \otimes \Cal M^t) = p^\ast \Cal M^t \otimes {\Cal O}_{\tilde X} (-E)$ is very ample, again for $t$ large enough.  So we can take $\Cal L$ as the linear system corresponding to the global sections of $\Cal I \otimes \Cal M^t$.

When $X$ is quasiprojective, we can apply the previous argument to its projective closure and restrict the obtained linear system to $X$. \qed
\enddemo
\bigskip
\proclaim
{5.3. Proposition}  Let $X$ be a quasiprojective normal variety and $p : \tilde X \rightarrow X$ the blowing--up of $X$ in $X_{\sing}$. Let $\pi_1$ and $\pi_2$ be two log resolutions of $X$ that factorize through $p$. Then there exists an effective Cartier divisor $D$ on $X$ with $X_{\sing} \subset \supp D$, such that for $i = 1,2$ the strict transform of $D$ by $\pi_i$ has normal crossings with the exceptional divisor of $\pi_i$, i.e. $\pi_i$ is also a log resolution of the pair $(X,D)$.
\endproclaim
\bigskip
\demo{Proof}  Say $\pi_i$ factorizes as $\pi_i : Y_i \overset h_i \to \longrightarrow \tilde X \overset p \to \longrightarrow X$.  Consider the linear system 
$\Cal L$ of Lemma 5.2; the induced system $h^\ast_i (p^\ast \Cal L - E)$ on $Y_i$ is also base--point free.  Take now $D$ as a general member of $\Cal L$.  Its strict transform by $\pi_i$ is a general member of the base--point free linear system $h^\ast_i (p^\ast \Cal L - E)$, and has thus  normal crossings with the exceptional divisor of $\pi_i$ by Bertini's Theorem (see e.g. [Jo, Theorem 6.10]). \qed \enddemo
\bigskip
\proclaim
{5.4. Theorem} Let $X$ be a quasiprojective $\Bbb Q$--Gorenstein variety.  Let $\pi : Y \rightarrow X$ be any log resolution of $X$ that factorizes through the blowing--up $p : \tilde X \rightarrow X$ of $X$ in $X_{\sing}$, and, using the notation of (0.1), such that all log discrepancies $a_i$ of exceptional divisors $E_i$ of $\pi$ with respect to $X$ are nonzero. Then
$$\Cal E(X) := \sum_{I \subset T} [E^\circ_I]\prod_{i \in I} \frac{L-1}{L^{a_i}-1}$$ does not depend on the chosen such resolution.
\endproclaim
\bigskip
\demo{Proof}  Let $\pi^\prime : Y^\prime \rightarrow X$ be another such resolution for which $E^\prime_i, i \in T^\prime,$ are the irreducible components of the exceptional divisor with log discrepancies $a^\prime_i$, and denote $E^{\prime \circ}_I := (\cap_{i \in I} E^\prime_i) \setminus (\cup_{\ell \not \in I} E^\prime_\ell)$ for $I \subset T^\prime$.  Take an effective Cartier divisor 
$D$ on $X$ as in Proposition 5.3, associated to $\pi$ and $\pi^\prime$.  Let $E_i, i \in T_s,$ denote the irreducible components of the strict transform of $D$ by $\pi$, and say $\pi^\ast D = \sum_{i \in T \cup T_s} N_i E_i$.

We consider the zeta function $\Cal Z(D,s)$ on the $\Bbb Q$--Gorenstein variety $X$ of [Ve2, \S6], associated to the effective Cartier divisor $D$, whose support contains $X_{\sing}$ as required there.  By the formula of [Ve2, \S6] for $\Cal Z(D,s)$ in terms of $\pi$  we have
$$L^{\dim X}\Cal Z(D,s) = \sum_{I \subset T \cup T_s} [E^\circ_I] \prod_{i \in I} \frac{L-1}{L^{a_i+sN_i}-1}\, .$$
Here the log discrepancies $a_i, i  \in T_s,$ are just $1$, and the notation $E^\circ_I$ should be clear. 

Now $\Cal Z(D,s)$ is an invariant of $X$, hence so is its specialization $\Cal Z(D,0)$.  Note that this specialization makes sense since all $a_i, i \in T \cup T_s$, are nonzero.  Clearly
$$L^{\dim X}\Cal Z(D,0) = \sum_{I \subset T \cup T_s} [E^\circ_I] \prod_{i \in I} \frac{L-1}{L^{a_i}-1} = \sum_{I \subset T}[E^\circ_I] \prod_{i \in I} \frac{L-1}{L^{a_i}-1} \tag $*$ $$     
because $a_i = 1$ for $i \in T_s$.  Analogously we obtain that
$$L^{\dim X}\Cal Z (D,0) = \sum_{I \subset T^\prime} [E^{\prime \circ}_I] \prod_{i \in I} \frac{L-1}{L^{a_i^\prime}-1}, \tag $*'$ $$ 
and so the right hand sides of ($*$) and ($*'$) are indeed equal. \qed
\enddemo
\bigskip
\noindent
{\bf 5.4.1.}  Restricting to log resolutions $\pi : Y \rightarrow X$ that factorize through the blowing--up of $X$ in $X_{\sing}$, we can say that $\Cal E (X)$ as above is a {\it (partial) stringy $\Cal E$--invariant} of those $X$ for which there exists such a resolution $\pi$ with all log discrepancies $a_i \ne 0$.
\bigskip
\noindent
{\bf 5.5.} Theorem 5.4 can be generalized to pairs $(X,B)$, where $X$ is a quasiprojective normal variety, $B = \sum_i b_i B_i$ is a $\Bbb Q$--divisor on $X$ where the $B_i$ are distinct and irreducible and all $b_i < 1$, and $K_X + B$ is $\Bbb Q$--Cartier.  Here we consider log resolutions $\pi : Y \rightarrow X$ of $(X,B)$ that factorize through the blowing--up of $X$ in $X_{\sing} \cup B_{\nnc}$, where $B_{\nnc}$ is the part of $\supp B$ in which $B$ is not a normal crossings divisor (this is the natural generalization of $X_{\sing}$ to pairs).  Then for those $(X,B)$ that admit such a log resolution $\pi$ for which all log discrepancies $a_i$ (with respect to $(X,B)$) are nonzero, we have that
$$\Cal E(X,B) := \sum_{I \subset T} [E^\circ_I]\prod_{i \in I} \frac{L-1}{L^{a_i}-1}$$
is an invariant of $(X,B)$, still using the notation of (0.1).

The proof is analogous and follows from the appropriate generalization of Proposition 5.3 and [Ve2, \S6] to pairs $(X,B)$.
\bigskip
\noindent
{\bf 5.6.} In the context of this section, but also more generally, the following question is important.
\bigskip
\noindent
{\it Problem.} Let $X$ be a (normal) variety and $\pi_1 : Y_1 \rightarrow X$ and $\pi_2 : Y_2 \rightarrow X$ two different log resolutions of $X$.  Does there exist an effective Cartier divisor $D$ on $X$ with $X_{\sing} \subset \supp D$, such that $\pi_1$ and $\pi_2$ are also log resolutions of $(X,D)$? (And analogously for a given pair $(X,B)$.)
\bigskip
For instance in [BL2, Remark 3.11], Borisov and Libgober are confronted with the same problem when they want to associate an elliptic genus to pairs $(X,B)$, with $X$ projective, which are not klt.

\enddocument